\documentclass[12pt]{amsart}
\usepackage{amssymb,amsmath,amscd}
\usepackage{latexsym,bm,bbm,mathrsfs}
\usepackage{hyperref,graphicx}
\usepackage{mathtools}
\usepackage{stmaryrd}
\usepackage[all,pdf]{xy}
\usepackage{graphicx}
\usepackage{color}
\graphicspath{ {./} }

\newtheorem{thm}{Theorem}[section]

\newtheorem{cor}[thm]{Corollary}
\newtheorem{prop}[thm]{Proposition}
\newtheorem{lemma}[thm]{Lemma}

\newtheorem{remark}[thm]{Remark}

\newtheorem{defn}[thm]{Definition}

\newcommand{\bbq}{{\mathbb{Q}}}

\newcommand{\bbc}{\mathbb{C}}
\newcommand{\bbf}{\mathbb{F}}
\newcommand{\bbh}{\mathbb{H}}

\newcommand{\Tr}{\operatorname{Tr}}

\newcommand{\ord}{\operatorname{ord}}

\title{Counting rational curves on K3 surfaces with finite group actions}

\author{Sailun Zhan}
\address{Department of Mathematics, Indiana University, Bloomington, IN, 47405, U.S.A.}
\email{zhans@indiana.edu}
\begin{document}

\begin{abstract}
G\"ottsche gave a formula for the dimension of the cohomology of Hilbert schemes of points on a smooth projective surface $S$. When $S$ admits an action by a finite group $G$, we describe the action of $G$ on the Hodge structure. In the case that $S$ is a K3 surface, each element of $G$ gives a trace on $\sum_{n=0}^{\infty}\sum_{i=0}^{\infty}(-1)^{i}H^{i}(S^{[n]},\mathbb{C})q^{n}$. When $G$ acts faithfully and symplectically on $S$, the resulting generating function is of the form $q/f(q)$, where $f(q)$ is a cusp form. We relate the Hodge structure of Hilbert schemes of points to the Hodge structure of the compactified Jacobian of the tautological family of curves over an integral linear system on a K3 surface as $G$-representations. Finally, we give a sufficient condition for a $G$-orbit of curves with nodal singularities not to contribute to the representation.
\end{abstract}

\maketitle

\section{Introduction}

Let $S$ be a smooth projective K3 surface over $\mathbb{C}$. In \cite{YZ96} and \cite{Bea99}, the number of rational curves in an integral linear system on $S$ is calculated using the relative compactified Jacobian. The idea is that the Euler characteristic of the relative compactified jacobian equals the number of maximally degenerate fibers if all rational fibers are nodal, and these are the rational curves we want. But the relative compactified jacobian is birational to the Hilbert scheme of points of $S$, the Euler characteristic of which is computed in \cite{Go90}. Hence we get the number of rational curves and the generating series:
\[
\sum_{n=0}^{\infty}N(n) t^{n}=\sum_{n=0}^{\infty}e({\overline{J^{n}(\mathcal{C}_{n})}})t^{n}=\sum_{n=0}^{\infty}e(S^{[n]}) t^{n}=\prod_{n=1}^{\infty}(1-t^{n})^{-24}=\frac{t}{\Delta(t)}
\] 
where  $N(n)$ is the number of rational curves contained in an $n$-dimensional linear system $|\mathcal{L}|$, $\mathcal{C}_{n}$ is the tautological family of curves over $|\mathcal{L}|$ with fibers being integral, $S^{[n]}$ is the Hilbert scheme of $n$ points of $S$ and $\Delta(t)=t\prod_{n\geq 1}(1-t^n)^{24}$ is the unique cusp form of weight 12 for $SL_{2}(\mathbb{Z})$.

In this paper $G$ will always be a finite group. We will consider a smooth projective K3 surface $S$ over $\mathbb{C}$ with a $G$-action, and ask whether we can prove a similar equality for $G$-representations. There are many different methods to obtain the generating series of the Euler characteristic of the Hilbert scheme of points of $S$. We follow the approach by Nakajima and Grojnowski \cite{Gr96}\cite{N97}\cite{N99}, which describe the sum of the cohomology groups of all the Hilbert schemes of points $\oplus_{n=0}^{\infty} H^{*}(S^{[n]})$ as the Fock space $\bbf(H^{*}(S))$ (See Section 2). We will consider the G-equivariant Hodge-Deligne polynomial for a smooth projective variety $X$
\[
E(X;u,v)=\sum_{p,q}(-1)^{p+q}[H^{p,q}(X,\bbc)]u^{p}v^{q},
\]
where the coefficients lie in the ring of virtual G-representations $R_{\bbc}(G)$, of which the elements are the formal differences of isomorphism classes of finite dimensional $\bbc$-representations of $G$. The addition is given by direct sum and the multiplication is given by tensor product. Define 
\[
F(X;u,v):=E(X;-u,-v)=\sum_{p,q}[H^{p,q}(X,\bbc)]u^{p}v^{q}.
\]

We may abbreviate $F(X;u,v)$ (resp. $E(X;u,v)$) as $F(X)$ (resp. $E(X)$).

The main results are as follows.

\begin{thm}\label{thm1}
Let $S$ be a smooth projective surface over $\bbc$ with a $G$-action. Let $S^{[n]}$ be the Hilbert scheme of n points of $S$. Then we have the following equality as virtual $G$-representations.
\[
\sum_{n=0}^{\infty}F(S^{[n]})t^{n}=\prod_{m=1}^{\infty}\prod_{p,q}\left({\sum_{i=0}^{h_{p,q}}(-1)^{i(p+q+1)}[\wedge^{i}H^{p,q}(S,\bbc)]u^{i(p+m-1)}v^{i(q+m-1)}t^{mi}} \right)^{(-1)^{p+q+1}},
\]
where $h_{p,q}$ are the dimensions of the Hodge pieces $H^{p,q}(S,\bbc)$.
\end{thm}

\begin{cor}\label{stable}
Let $S$ be a smooth projective surface over $\mathbb{C}$ with a $G$-action. If we fix $p,q\geq 0$, then $H^{p,q}(S^{[n]},\mathbb{C})$ become stable for $n\geq p+q$ as $G$-representations.
\end{cor}

\begin{defn}
Let $X$ and $Y$ be smooth projective varieties over $\bbc$ with G-actions. They are called \emph{G-equivariant K-equivalent} if there exists a smooth projective variety $Z$ with a G-action and G-equivariant birational morphisms $f\colon Z\to X$ and $g\colon Z\to Y$ such that $f^{*}\omega_{X}\cong g^{*}\omega_{Y}.$
\end{defn}

\begin{thm}\label{K-equiv}
Let $X$ and $Y$ be smooth projective varieties over $\mathbb{C}$ with $G$-actions. If $X$ and $Y$ are G-equivariant K-equivalent, then
\[
H^{*}(X,\bbq)\cong H^{*}(Y,\bbq)
\]
as Hodge structures with G-actions.
\end{thm}

Recall that a linear system $|\mathcal{L}|$ is called an integral linear system if every effective divisor in it is integral. $|\mathcal{L}|$ is called $G$-stable if $G$ induces an action on the projective space $|\mathcal{L}|$, which means $G$ maps an effective divisor in $|\mathcal{L}|$ again to an effective divisor in $|\mathcal{L}|$.

\begin{cor}\label{euler}
 Let $S$ be a smooth projective K3 or abelian surface over $\mathbb{C}$ with a $G$-action, and let $\mathcal{C}_n$ be the tautological family of curves over any $n$-dimensional integral $G$-stable linear system. Then we have the following equality as virtual $G$-representations, where $\overline{J^{n}(\mathcal{C}_n)}$ is the relative compactified jacobian.
\[
\sum_{n=0}^{\infty}E(\overline{J^{n}(\mathcal{C}_n)})t^{n}=\prod_{m=1}^{\infty}\prod_{p,q}\left({\sum_{i=0}^{h_{p,q}}(-1)^{i}[\wedge^{i}H^{p,q}(X,\bbc)]u^{i(p+m-1)}v^{i(q+m-1)}t^{mi}} \right)^{(-1)^{p+q+1}}
\]
\end{cor}

\begin{remark}
Note that we are fixing the surface $S$ here, so the equality above should be understood as: if $S$ admits an $n$-dimensional integral $G$-stable linear system, then $E(\overline{J^{n}(\mathcal{C}_n)})$ equals the coefficient of $t^n$ on the right hand side. 
\end{remark}

Recall that for a complex K3 surface $S$ with an automorphism $g$ of finite order $n$, $H^{0}(S,K_{S})=\mathbb{C}\omega_{S}$ has dimension 1, and we say $g$ acts symplectically on $S$ if it acts trivially on $\omega_{S}$, and $g$ acts non-symplectically otherwise, namely, $g$ sends $\omega_{S}$ to $\zeta_{n}^{k}\omega_{S}$, $0<k<n$, where $\zeta_{n}$ is a primitive $n$-th root of unity.

We denote by $[e(X)]$ the alternating sum of cohomology groups $E(X;1,1)\in R_{\bbc}(G)$. For the right hand side of the equality in Corollary \ref{euler}, we have

\begin{thm}\label{symplectic}
Let $G$ be a finite group which acts faithfully and symplectically on a complex K3 surface $S$. Then
\[
\sum_{n=0}^{\infty}\Tr(g,[e(S^{[n]})])t^{n}=\exp\left(\sum_{m=1}^{\infty}\sum_{k=1}^{\infty}\frac{\epsilon(\ord(g^{k}))t^{mk}}{k}\right)
\]
for all $g\in G$, where $\epsilon(n)=24\left(n\prod_{p|n}\left(1+\frac{1}{p}\right)\right)^{-1}$. In particular, if $G$ is generated by a single element $g$ of order $N\leq 8$, then we have
\begin{center}
\begin{tabular}{c||c}
$N$ & $\sum_{n=0}^{\infty}\Tr(g,[e(S^{[n]})])t^{n}$ \\
\hline\hline
$1$ & $t/\eta(t)^{24}$ \\
$2$ & $t/\eta(t)^{8}\eta(t^{2})^{8}$ \\
$3$ & $t/\eta(t)^{6}\eta(t^{3})^{6}$ \\
$4$ & $t/\eta(t)^{4}\eta(t^{2})^{2}\eta(t^{4})^{4}$ \\
$5$ & $t/\eta(t)^{4}\eta(t^{5})^{4}$ \\
$6$ & $t/\eta(t)^{2}\eta(t^{2})^{2}\eta(t^{3})^{2}\eta(t^{6})^{2}$ \\
$7$ & $t/\eta(t)^{3}\eta(t^{7})^{3}$ \\
$8$ & $t/\eta(t)^{2}\eta(t^{2})\eta(t^{4})\eta(t^{8})^{2}$ 
\end{tabular}
\end{center}
where $\eta(t)=t^{1/24}\prod_{n=1}^{\infty}(1-t^{n})$.
\end{thm}

\begin{remark}
Notice that if $g$ acts symplectically on $S$, then $g$ has order $\leq 8$ by \cite[Corollary 15.1.8]{H16}. We know that the generating series for the {topological Euler characteristic} $\sum_{n=0}^{\infty}e(S^{[n]})t^{n}=\sum_{n=0}^{\infty}\Tr(1,[e(S^{[n]})])t^{n}$ equals $t/\Delta(t)$, where $\Delta(t)=\eta(t)^{24}$ is a level 1 cusp form of weight 12. But by Theorem 1.6, we deduce that when an element $g$ of order $N$ acts faithfully and symplectically on $S$, we have $\sum_{n=0}^{\infty}\Tr(g,[e(S^{[n]})])t^{n}=t/F(t)$, where $F(t)$ is also a cusp form for $\Gamma_{0}(N)$ of weight $\lceil\frac{24}{N+1}\rceil$ (\cite[Proposition 5.9.2]{CS17}). This coincides with the results by Jim Bryan and \'Ad\'am Gyenge in \cite{BG19} when $G$ is a cyclic group. See also \cite[Lemma 3.1]{BO18}.
\end{remark}

\begin{thm}\label{nonsymplectic}
Let $G=\left<g\right>$ be a finite group generated by an automorphism $g$ of order $p$, where $p$ is a prime number. Suppose $g$ acts non-symplectically on a complex K3 surface $S$. Then we have
\[
\sum_{n=0}^{\infty}\Tr(1,[e(S^{[n]})])t^{n}=\left(\prod_{m=1}^{\infty}(1-t^{m})\right)^{-24}
\]
\[
\sum_{n=0}^{\infty}\Tr(g,[e(S^{[n]})])t^{n}=\left(\prod_{m=1}^{\infty}(1-t^{m})\right)^{dp-24}\left(\prod_{m=1}^{\infty}(1-t^{mp})\right)^{-d}
\]
for all $g\in G$, $g\neq 1$, where $d=\frac{{\rm{rank}}\,T(g)}{p-1}$, and $T(g):=(H^{2}(S,\mathbb{Z})^{g})^{\bot}$ is the orthogonal complement of the $g$-invariant sublattice.
\end{thm}

For the left hand side of the equality in Corollary \ref{euler}, recall that in \cite{Bea99}, for the curve $C$ in $\mathcal{C}_n$, the {topological Euler characteristic} $e(\overline{J^{n}(C)})=0$ if the normalization $\tilde{C}$ has genus $\geq$ 1, and the {topological Euler characteristic} $e(\overline{J^{n}(C)})=1$ if $C$ is a nodal rational curve. Hence intuitively that is why $e(\overline{J^{n}(\mathcal{C}_{n})})$ counts the number of rational curves in $\mathcal{C}_{n}$ if we assume all rational curves in $\mathcal{C}_{n}$ are nodal. But in our situation, $e(\overline{J^{n}(C)})=0$ does not mean $[e(\overline{J^{n}(C)})]=0$ as $G$-representations. Hence non-rational curves may also contribute to $[e({\overline{J^{n}(\mathcal{C}_n)}})]$, and certain $G$-orbits of curves contribute certain representations (See Example 1, 2 and 3). Nevertheless, we show that a $G$-orbit of curves with nodal singularities will contribute nothing if the normalization of the curve quotient by its stablizer is not rational. By this method, we are able to understand certain $G$-orbits in the linear system. We denote by $[e(X)]$ the alternating sum of the compactly supported $l$-adic cohomology $\sum_{n=0}(-1)^{n}[H^{n}_{c}(X_{\overline{\bbf}_{p}},\mathbb{Q}_{l})]$ when we are in the situation of characteristic $p$.

\begin{thm}\label{curve1}
Let $C$ be an integral curve over $\overline{\bbf}_{p}$ with nodal singularities and a $G$-action. Suppose $p\nmid|G|$. Denote by $\tilde{C}$ its normalization. If $\tilde{C}/\left<g\right>$ is not a rational curve {\rm(}i.e. $\mathbb{P}^{1}${\rm)} over $\overline{\bbf}_{p}$ for every $g\in G$, then $[e(\overline{J^{n}C})]=0$ as $G$-representations.
\end{thm}

\begin{cor}\label{curve2}
Let $C$ be an integral curve over $\overline{\bbf}_{p}$ with nodal singularities and a $G$-action. Suppose $p\nmid|G|$. Denote by $\tilde{C}$ its normalization. If $\tilde{C}/G$ is not a rational curve {\rm(}i.e. $\mathbb{P}^{1}${\rm)} over $\overline{\bbf}_{p}$, then $[e(\overline{J^{n}C})]=0$ as $G$-representations.
\end{cor}

By the above discussions, we show that the representation $[e(\overline{J^{n}(\mathcal{C}_{n})})]$ actually `counts' the curves in $\mathcal{C}_{n}$ whose normalization quotient by its stablizer is rational (See Example 1, 2, and 3).

This paper is organized as follows. In Section 2, we recall the Nakajima operators on the cohomology groups of Hilbert schemes of points, and we show that the theory works in the G-equivariant settings. In Section 3, we work with the G-equivariant Grothendieck ring and prove Theorem \ref{K-equiv} via motivic integration. In Section 4, we deal with compactified jacobians and prove Theorem \ref{curve1} and Corollary \ref{curve2}. In Section 5, we prove Corollary \ref{euler}, Theorem \ref{symplectic} and Theorem \ref{nonsymplectic} by using the results in previous sections. Then we give three explicit examples when $G$ equals $\mathbb{Z}/2\mathbb{Z}$ or $\mathbb{Z}/3\mathbb{Z}$. Finally, when $G=PSL(2,7), A_{6}, A_{5}\ {\rm or}\ S_{5}$, we show that the smooth projective curve $C$ over $\mathbb{C}$ with a faithful $G$-action must be of specific kind if there exists $g\in G$ such that $C/{\left<g\right>}=\mathbb{P}^{1}$. 

\section{Hilbert schemes of points}

Let $X^{[n]}$ denote the component of the Hilbert scheme of a projective scheme $X$ parametrizing subschemes of length $n$ of $X$. For properties of Hilbert scheme of points, see references \cite{I77}, \cite{Go94} and \cite{N99}. The following theorem is proved for smooth projective surfaces over $\mathbb{C}$ in \cite{Go90}, and for smooth quasi-projective surfaces over $\mathbb{C}$ in \cite{GS93}.

\begin{thm}
Let $S$ be a smooth quasi-projective surface over $\bbc$. Then the generating function of the Poincar\'e polynomials of the Hilbert scheme $S^{[n]}$ is given by 
\[
\sum_{n=0}^{\infty} P(S^{[n]},z)t^n=\prod_{m=1}^{\infty} \frac{(1+z^{2m-1}t^m)^{b_{1}(S)}(1+z^{2m+1}t^{m})^{b_{3}(S)}}{(1-z^{2m-2}t^m)^{b_{0}(S)}(1-z^{2m}t^m)^{b_{2}(S)}(1-z^{2m+2}t^m)^{b_{4}(S)}}.
\]
\end{thm}

One can also prove the above theorem for a smooth projective surface $S$ by constructing an action of the Heisenberg algebra on the direct sum of the cohomology groups of all the $S^{[n]}$ \cite{Gr96}\cite{N97}\cite{N99}. We recall some of the constructions following \cite{EG00}.

A vector space $V$ over $\bbq$ is called a super vector space if there is a decomposition $V=V^+\oplus V^-$ of $V$ into an even and an odd part. For any super vector space $V$, one can construct an algebra $\bbf(V)$, which is called the \emph{Fock sapce}. There is an isomorphism of graded vector space
\[
\bbf(V)\cong\bigotimes_{m=0}^{\infty}S(V^{+}\otimes t^m)\otimes\wedge(V^{-}\otimes t^{m}),
\]
where 
\[
S(V):=\bigoplus_{i\geq 0} S^{i}(V),\ \ \ \ \wedge(V):=\bigoplus_{i\geq 0}\wedge^{i}(V)
\]
are the symmetric and alternating algebra on $V$. The grade is given by the exponent of the $t$ powers. The element in the infinite tensor product should be understood as a finite linear combinition of the elements in some finite tensor product.

One can also construct the \emph{current algebra (Heisenberg algebra)} $\mathbb{S}(V)$, which acts on the Fock space $\bbf(V)$. For a smooth projective surface $S$ over $\bbc$, the cohomology groups $H^{*}(S)=H^{+}(S,\bbq)\oplus H^{-}(S,\bbq)$ is a super vector space according to the degree of the cohomology. Hence we can construct the corresponding Fock space $\bbf(H^{*}(S))$ and the current algebra $\mathbb{S}(H^{*}(S))$. Define
\[
\bbh^{+}(S):=\bigoplus_{n\geq 0, i\geq 0} H^{2i}(S^{[n]},\bbq),
\]
\[
\bbh^{-}(S):=\bigoplus_{n\geq 0, i\geq 0} H^{2i+1}(S^{[n]},\bbq),
\]
and let $\bbh(S):=\bbh^{+}(S)\oplus\bbh^{-}(S)$.

It turns out that $\bbh(S)$ is a $\mathbb{S}(H^{*}(S))$-module and the following theorem describes the cohomology groups of Hilbert schemes of points in terms of the cohomology groups of $S$.

\begin{thm}\cite[Theorem 5.6]{EG00}\label{Nak}
The Fock space $\bbf(H^{*}(S))$ and the space $\bbh(S)$ are isomorphic as $\mathbb{S}(H^{*}(S))$-modules.
\end{thm}

The map from $\bbf(H^{*}(S))$ to $\bbh(S)$ is given by sending $u\otimes t^{i}$ to $q_{i}[u]\bf{1}$, where $q_{i}[u]$ are the Nakajima creation operators, which are defined in \cite[\S 4]{EG00}. We observe that the definition of $q_{i}[u]$ is geometric and only involves pull back/push forward operators on cohomology/homology groups and the cap product/Poincar\'e duality homomorphism. But all of those operators are G-equivariant and are morphisms of mixed Hodge structures because of the functoriality of mixed Hodge structures and \cite[Proposition 6.18]{PS08}. In particular, the isomorphism in Theorem \ref{Nak} is G-equivariant, and we obtain 
\begin{equation}\label{creation}
q_{i}[u]\alpha\in H^{p+r+i-1,q+s+i-1}(X^{[n+i]})
\end{equation}
for $i\geq 1$, $u\in H^{r,s}(X)$ and $\alpha\in H^{p,q}(X^{[n]})$.

\begin{proof}[Proof of Theorem \ref{thm1}]
Since a G-representation is determined by the trace, we can assume that G is a cyclic group without loss of generality. We deduce that
\[
\begin{aligned}
\bbf(H^{*}(S))&=\bigotimes_{m=0}^{\infty}S(H^{+}(S)\otimes t^m)\otimes\wedge(H^{-}(S)\otimes t^{m})\\
&=\prod_{m=1}^{\infty}\prod_{p,q}\prod_{i=1}^{h_{p,q}}\left(1+(-1)^{p+q+1}H_{i}^{p,q}(S)\otimes t^{m}\right)^{(-1)^{p+q+1}},
\end{aligned}
\]
as graded G-representation, where $H^{p,q}_{i}(S)$ are eigenspaces of $G$ acting on $H^{p,q}(S)$. Combining the fact that $\bbf(H^{*}(S))$ is mapped G-equivariantly to $\bbh(S)$ as $\mathbb{S}(H^{*}(S))$-modules and the relation \refeq{creation}, we deduce that
\[
\begin{aligned}
\sum_{n=0}^{\infty}F(S^{[n]})t^{n}&=\prod_{m=1}^{\infty}\prod_{p,q}\prod_{i=1}^{h_{p,q}}\left(1+(-1)^{p+q+1}[H^{p,q}_{i}(S)]u^{p+m-1}v^{q+m-1}t^{m}\right)^{(-1)^{p+q+1}}\\
&=\prod_{m=1}^{\infty}\prod_{p,q}\left({\sum_{i=0}^{h_{p,q}}(-1)^{i(p+q+1)}[\wedge^{i}H^{p,q}(S,\bbc)]u^{i(p+m-1)}v^{i(q+m-1)}t^{mi}} \right)^{(-1)^{p+q+1}}.
\end{aligned}
\]
as virtual G-representations.
\end{proof}

If we only consider the cohomology groups, then Theorem \ref{thm1} implies the following statements.
\begin{cor}\label{trace}
Let $S$ be a smooth projective surface over $\bbc$. Then
\[
\sum_{n=0}^{\infty}[e(S^{[n]})]t^{n}=\prod_{m=1}^{\infty}\prod_{k=0}^{4}\left({\sum_{i=0}^{b_{k}}(-1)^{i}[\wedge^{i}H^{k}(S,\bbc)]t^{mi}} \right)^{(-1)^{k+1}}
\]
as virtual G-representations, where $b_{k}$ are the Betti numbers of $S$, and
\[
\sum_{n=0}^{\infty}\Tr(g,[e(S^{[n]})])t^{n}=\prod_{m=1}^{\infty}\left(\frac{(\prod_{i=1}^{b_1}(1-{g_{1,i}}t^{m}))(\prod_{i=1}^{b_3}(1-{g_{3,i}}t^{m}))}{(1-t^{m})(\prod_{i=1}^{b_2}(1-{g_{2,i}}t^{m}))(1-t^{m})}\right)
\] 
\[
=\exp\left(\sum_{m=1}^{\infty}\sum_{k=1}^{\infty}\frac{t^{mk}}{k}\left(1-\sum_{i=1}^{b_{1}}g_{1,i}^{k}+\sum_{i=1}^{b_{2}}g_{2,i}^{k}-\sum_{i=1}^{b_{3}}g_{3,i}^{k}+1\right)\right),
\]
where $g_{j,i}$ are the eigenvalues of $g$ acting on $H^{j}(S,\bbc)$, $0\leq j\leq 4$. 
\end{cor}

We will need the above expression in Section 5. Now we prove Corollary \ref{stable}.

\begin{proof}[Proof of Corollary \ref{stable}]
We let
\[
G(u,v,t):=(1-t)\sum_{n=0}^{\infty}F(S^{[n]})t^{n}
\]
\[
=(1-t)\prod_{m=1}^{\infty}\prod_{p,q}\left({\sum_{i=0}^{h_{p,q}}(-1)^{i(p+q+1)}[\wedge^{i}H^{p,q}(S,\bbc)]u^{i(p+m-1)}v^{i(q+m-1)}t^{mi}} \right)^{(-1)^{p+q+1}}
\]
We denote by $a_{i,j,k}$ the coefficient of $u^{i}v^{j}t^{k}$. If $i+j<k$, then $a_{i,j,k}(G(u,v,t))=0$ as $G$-representation. Now fix $p,q$ and take $n\geq p+q$. Then we deduce that
\[
\begin{aligned}
H^{p,q}(S^{[n]},\mathbb{C})&=a_{p,q,n}\Bigl(\Bigl(\sum_{k=0}^{\infty}t^{k}\Bigr)G(u,v,t)\Bigr)\\
&=\sum_{k=0}^{n}a_{p,q,k}(G(u,v,t))\\
&=\sum_{k=0}^{\infty}a_{p,q,k}(G(u,v,t))\\
&=a_{p,q,0}(G(u,v,1)).
\end{aligned}
\]
Notice that $a_{p,q,0}(G(u,v,1))$ is a representation independent of $n$. Hence $H^{p,q}(S^{[n]},\mathbb{C})$ become stable as $G$-representations for $n\geq p+q$. 
\end{proof}

\section{K-equivalent smooth projective varieties}

Recall that two smooth projective varieties $X$ and $Y$ are called K-equivalent if there exists a smooth projective variety $Z$  and birational morphisms $f\colon Z\to X$ and $g\colon Z\to Y$ such that $f^{*}\omega_{X}\cong g^{*}\omega_{Y}.$

The following result is proved by Kontsevich.
\begin{thm}
Let $X$ and $Y$ be smooth projective complex K-equivalent varieties. Then
\[
h^{p,q}(X)=h^{p,q}(Y)
\]
for all $p,q$. 
\end{thm}
The proof uses motivic integrations and actually shows that $[X]=[Y]\in\widehat{\mathcal{M}_{\bbc}}$ in a localized and completed Grothendieck ring of varieties. For the proof in the non-equivariant case, see \cite{Bl11} and \cite{L09}. We will prove a G-equivariant version of the theorem.

\begin{defn}
The G-equivariant Grothendieck ring of varieties $K_{0}^{G}(\emph{Var}_{\bbc})$ is the quotient of the free abelian group generated by isomorphism classes $[X]$ of G-varieties with good G-actions (i.e. every orbit is contained in an affine open subset) by the following relations:

(1) $[X]=[Y]+[X\backslash Y]$, where $Y$ is a G-equivariant closed subvariety of $X$.

(2) $[X\times\mathbb{A}^{n}_{\bbc}]=[X][\mathbb{A}^{n}_{\bbc}]$ if the G-action on $[X\times\mathbb{A}^{n}_{\bbc}]$ lifts the G-action on $[X]$, and the G-action on $[\mathbb{A}^{n}_{\bbc}]$ is trivial.
\end{defn}

This definition was used in \cite[\S 2.1]{LN20} and \cite[2.9.]{DL02}. Any G-action on a quasi-projective variety is a good G-action. For the equivariant motivic integration, see the discussions in \cite{LN20}, which follows \cite{DL99}.

Let $\mathcal{M}_{\bbc}^G:=K_{0}^{G}(\text{Var}_{\bbc})[\mathbb{L}^{-1}]$ be the localization of $K_{0}^{G}(\text{Var}_{\bbc})$ in the class $\mathbb{L}$ with a trivial G-action. Let $\widehat{\mathcal{M}_{\bbc}^{G}}$ be the completion with respect to the virtual dimension filtration. Denote the image of ${\mathcal{M}_{\bbc}^{G}}$ in $\widehat{\mathcal{M}_{\bbc}^{G}}$ by $\overline{\mathcal{M}_{\bbc}^{G}}$.

For a smooth projective variety $X$, we denote the $m$-th jet space of $X$ by $J_{m}(X)$ and the arc space of $X$ by $J_{\infty}(X)$. 

We will need an equivariant transformation rule, which is essentially proved in \cite[Theorem 3.1]{LN20}. If we let the measure take values in $\widehat{\mathcal{M}_{\bbc}^{G}}$, then we deduce the following.

\begin{thm}\label{trans}
Let $h\colon X\to Y$ be a G-equivariant proper birational morphism between smooth complex varieties with good G-actions. If $F:J_{\infty}(Y)\to\mathbb{N}\cup\{\infty\}$ is a G-invariant simple function, then
\[
\int_{J_{\infty}(Y)}\mathbb{L}^{-F}=\int_{J_{\infty}(X)}\mathbb{L}^{-(F\circ h_{\infty}+\ord_{h^{*}\omega_{Y}})},
\]
where the equality holds in $\widehat{\mathcal{M}_{\bbc}^{G}}$, and $h_{\infty}\colon J_{\infty}(X)\to J_{\infty}(Y)$ is the map induced by $h$.
\end{thm}

\begin{proof}[Proof of Theorem \ref{K-equiv}]
Without loss of generality, we can assume $G=\langle\mu\rangle$ is cyclic. By \cite[3.4]{DL02} (taking the Hodge realization of $K_{0}(\text{Mot}_{\bbc,E})$), we deduce that there is a map from $\overline{\mathcal{M}_{\bbc}^{G}}$ to the Grothendieck ring of Hodge structures with a $\langle\mu\rangle$-action. Hence it suffices to prove that
\[
[X]=[Y]\in\widehat{{\mathcal{M}}_{\bbc}^{G}}.
\]
Recall that if $X$ and $Y$ are G-equivariant K-equivalent, there exists a smooth projective variety $Z$ with a G-action and G-equivariant birational morphisms $f\colon Z\to X$ and $g\colon Z\to Y$ such that $f^{*}\omega_{X}\cong g^{*}\omega_{Y}.$ Apply Theorem \ref{trans} to the function $F\equiv 0$ with respect to both $f$ and $g$. Suppose $\text{dim}X=\text{dim}Y=n$. Then
\[
\frac{[X]}{\mathbb{L}^{n}}=\int_{J_{\infty}(X)}\mathbb{L}^{-0}=\int_{J_{\infty}(Z)}\mathbb{L}^{-\ord_{f^{*}\omega_{X}}}=\int_{J_{\infty}(Z)}\mathbb{L}^{-\ord_{g^{*}\omega_{Y}}}=\int_{J_{\infty}(Y)}\mathbb{L}^{-0}=\frac{[Y]}{\mathbb{L}^{n}},
\]
which implies
\[
[X]=[Y]\in\widehat{{\mathcal{M}}_{\bbc}^{G}}.
\]
\end{proof}

\begin{cor}\label{calabi-yau}
Let $X$ and $Y$ be smooth projective algebraic varieties over $\bbc$ with G-actions. If $X$ and $Y$ have trivial canonical bundles and there is a G-equivariant birational map $f\colon X\to Y$, then
\[
F(X;u,v)=F(Y;u,v).
\]
\end{cor}

\begin{proof}
Since there is an equivariant birational map between smooth projective varieties $X$ and $Y$, we deduce that $X$ and $Y$ are G-equivariant K-equivalent by the equivariant weak factorization theorem (\cite[remark 2 after Theorem 0.3.1]{AKMW02} or \cite[Theorem 5.2.1]{Ma00}). Recall that 
\[
F(X;u,v):=\sum_{p,q}[H^{p,q}(X,\bbc)]u^{p}v^{q}.
\]
Hence by Theorem \ref{K-equiv}, it follows that $F(X;u,v)=F(Y;u,v)$.

\end{proof}

\section{Compactified Jacobians }

We will need to consider varieties over finite fields in this section, since we do not know whether Lemma \ref{bundle} is true for singular cohomology in characteristic 0.

Recall some facts from \cite{AK76}, \cite{A04} and \cite{EGK00}. Let $C/S$ be a flat projective family of integral curves. By a torsion-free rank-1 sheaf $\mathcal I$ on $C/S$, we mean an $S$-flat coherent $\mathcal{O}_C$-module $\mathcal I$ such that, for each point $s$ of $S$, the fiber $\mathcal{I}_s$ is a torsion-free rank-1 sheaf on the fiber $C_s$. We say that $\mathcal I$ is of degree $n$ if $\chi(\mathcal{I}_{s})-\chi(\mathcal{O}_{C_{s}})=n$ for each $s$.

Given $n$, consider the \'etale sheaf associated to the presheaf that assigns to each locally Noetherian $S$-scheme $T$ the set of isomorphism classes of torsion-free rank-1 sheaves of degree $n$ on $C_{T}/T$. This sheaf is representable by a projective $S$-scheme, denoted $\bar{J}_{C/S}^{n}$. It contains $J^{n}:=\rm{Pic}^{n}_{C/S}$ as an open subscheme. For every $S$-scheme $T$, we have a natural isomorphism $\bar{J}_{C_{T}/T}^{n}=\bar{J}_{C/S}^{n}\times T$. If $S=\rm{Spec}$$k$ for an algebraically closed field $k$, we denote $\bar{J}_{C/S}^{n}$ by $\overline{J^{n}C}$.

Recall that at the beginning we are considering $\mathcal{C}$, which is the tautological family of curves over an $n$-dimensional integral $G$-stable linear system. Since $\mathcal{C}$ has a stratification according to the geometric genus of the fibers and the $G$-action (see $\S 6$), we can temporarily focus our attention on $\overline{J^{n}C}$ for a single singular curve $C$ with a $G$-action (note that our $G$-action on $\overline{J^{n}C}$ is given by pushing forward the torsion-free rank-1 sheaves). This is reasonable since we have

\begin{lemma}\label{strata}
Let $X$ be an algebraic variety over $\overline{\bbf}_{p}$ with a $G$-action, $U$ an open subvariety of $X$, $Z:=X\backslash U$ the closed subvariety. If both $U$ and $Z$ are $G$-stable, then  
\[
[e(X)]=[e(U)]+ [e(Z)]
\]
as virtual $G$-representations, i.e. in $R_{\mathbb{Q}_{l}}(G)$, and $[e(X)]:=\sum_{n=0}^{\infty}(-1)^{n}[H_{c}^{n}(X,\mathbb{Q}_{l})]$ is the alternating sum of the compactly supported $l$-adic cohomology groups.
\end{lemma}

\begin{proof}
One way to see this is to consider the bounded long exact sequence $0\rightarrow H^{0}_{c}(U)\rightarrow H^{0}_{c}(X)\rightarrow H^{0}_{c}(Z)\rightarrow ...$ and check that Tr$(g)=0$ on $-[e(U)]+[e(X)]-[e(Z)]$ for every $g$. Another way is to use the description of Tr$(g)$ on $e(X)$ without involving cohomology \cite[Appendix(h)]{Ca85}. 
\end{proof}

Now we have an integral curve $C$ over $\overline{\bbf}_{p}$. Recall that $\overline{J^{n}C}$ parametrizes the isomorphism classes of torsion-free rank-1 sheaves of degree $n$ on $C$, and we have the following facts \cite{Bea99}.

\begin{prop}\label{jacobian facts}
Let $C$ be an integral curve over an algebraically closed field $k$.

(1) If $L\in \overline{J^{n}C}$ is a non-invertible torsion-free rank 1 sheaf, then $L=f_{*}L'$, where $L'$ is some invertible sheaf on some partial normalization $f:C'\rightarrow C$.

(2) If $f:C'\rightarrow C$ is a partial normalization of $C$, then the morphism $f_{*}:\overline{J^{n}C'}\rightarrow\overline{J^{n}C}$ is a closed embedding.
\end{prop}

Using these two facts, we obtain the following corollary.

\begin{cor}\label{jac strata}
Let all the singularities of an integral curve $C$ be nodal singularities. Then $\overline{J^{n}C}$ has the following stratification
\[
\overline{J^{n}C}=\coprod_{C'\rightarrow C} {J^{n}C'}
\]
where $\overline{J^{n}C}$ parametrizes rank-1 torsion-free sheaves of degree n, and $C'$ goes through all partial normalizations of $C$ {\rm(}including $C$ itself{\rm)}.
\end{cor}

Now let $J^{n}C'$ be some stratum which is preserved by $G$. We want to calculate the $G$-representation $[e(J^{n}C')]$. Here we need to make use of the short exact sequence of algebraic groups
\[
0\rightarrow L\rightarrow J^{n}C'\rightarrow J^{n}\tilde{C'}\rightarrow 0,
\]
where $L$ is a smooth connected linear algebraic group \cite[$\S$9 Corollary 11]{BLR90}, and $\widetilde{C'}$ is the normalization of $C'$. Since $L$ is linear, we have that $J^{n}C'$ is a principal Zariski fiber bundle over $J^{n}\tilde{C'}$ \cite[Chapter VII, Proposition 6]{Se88}. Now we need to prove the following lemma, which is used to prove lemma \ref{bundle}.

\begin{lemma}\label{Weil}
Let $X$ and $Y$ be two smooth projective varieties over $\overline{\bbf}_{p}$ with finite group $G$-actions. Suppose $X,Y$ and the actions of G can be defined over $\bbf_{q}$, where $q$ is a $p$ power. If $|{X}(\overline{\mathbb{F}}_{p})^{gF_{q^n}}|=|{Y}(\overline{\mathbb{F}}_{p})^{gF_{q^n}}|$ for every $n\geq 1$ and $g\in G$, then $H^{i}(X,\bbq_{l})\cong H^{i}(Y,\bbq_{l})$ as $G$-representations for every $i\geq 0$.
\end{lemma}

\begin{proof}
Fix $g\in G$. Denote by $F_{q}$ the geometric Frobenius over $\bbf_{q}$. Since the finite group action is defined over $\bbf_{q}$, the action $g$ commutes with $F_{q}$ and the action of $g$ on the cohomology group is semisimple. There exists a basis of the cohomology group such that the actions of $g$ and $F_{q}$ are in Jordan normal forms simultaneously. Let $\alpha_{i,j},j=1,2,...,a_i$ (resp.  $\beta_{i,j},j=1,2,...,b_i$) denote the eigenvalues of $F_q$ acting on $H^{i}({X}, \mathbb{Q}_l)$ (resp. $H^{i}({Y}, \mathbb{Q}_l)$) in such a basis, where $a_i$ (resp. $b_i$) is the $i$-th betti number. Let $c_{i,j},j=1,2,...,a_i$ (resp.  $d_{i,j},j=1,2,...,b_i$) denote the eigenvalues of $g$ acting on the same basis of $H^{i}(X, \mathbb{Q}_l)$ (resp. $H^{i}({Y}, \mathbb{Q}_l)$). Then the Grothendieck trace formula (\cite[Prop.3.3]{DL76} and \cite[Appendix(h)]{Ca85}) implies that
\[
|{X}(\overline{\mathbb{F}}_{p})^{gF_{q^n}}|=\sum_{i=0}^{\infty}(-1)^{i}\text{Tr}((gF_{q^{n}})^{*},H^{i}({X},\mathbb{Q}_{l}))
\]

Since $|{X}(\overline{\mathbb{F}}_{p})^{gF_{q^n}}|=|{Y}(\overline{\mathbb{F}}_{p})^{gF_{q^n}}|$ for every $n\geq 1$, we have 
\[
\sum_{i=0}^{\infty}(-1)^{i}\sum_{j=1}^{a_i}c_{i,j}\alpha_{i,j}^{n}=\sum_{i=0}^{\infty}(-1)^{i}\sum_{j=1}^{b_i}d_{i,j}\beta_{i,j}^{n}
\]
for every $n\geq 1$. By linear independence of the characters $\chi_{\alpha}:\mathbb{Z}^{+}\rightarrow\mathbb{C}, n\mapsto \alpha^n$ and the fact that $\alpha_{i,j},\beta_{i,j},j=1,2,...$ all have absolute value $q^{i/2}$ by Weil's conjecture, we deduce that $a_{i}=b_{i}$ and $\sum_{j=1}^{a_i}c_{i,j}=\sum_{j=1}^{b_i}d_{i,j}$ for each $i$. But since $g$ is arbitrary, this implies that the $G$-representations $H^{i}({X},\mathbb{Q}_{l})$ and $H^{i}({Y},\mathbb{Q}_{l})$ are the same. 
\end{proof}
\begin{lemma}\label{bundle} 
Let $B$, $E$ and $F$ be separated schemes of finite type over $\overline{\bbf}_{p}$. Suppose $E$ is a Zariski-locally trivial fiber bundle over $B$ with fiber $F$ and let $G$ be a finite group acting on $E$ and $B$, the action of which is compatible with the projection $\pi:E\rightarrow B$. If $B$ is irreducible, then we have
\[
[e(E)]=[e(B)][e(F)]
\]
as virtual $G$-representations, i.e. in $R_{\mathbb{Q}_{l}}(G)$, where the action of $g\in G$ on $H_{c}^{*}(F,\mathbb{Q}_{l})$ is given by $H_{c}^{*}(\{g(b)\}\times F,\mathbb{Q}_{l})\rightarrow H_{c}^{*}(\{b\}\times F,\mathbb{Q}_{l})$ by choosing any closed point $b\in B$ in a $G$-invariant open subset which trivializes the bundle.
\end{lemma}

\begin{proof}
We first deal with the case when $E=B\times F$ is a trivial bundle. We begin with a homotopy argument. Fix $g\in G$. By assumption, we have a commutative diagram
\[
\begin{CD}
B\times F @>g>> B\times F\\
@V\pi VV @VV\pi V\\
B @>g>> B
\end{CD}
\]
Hence we have a map $\phi=(g,\pi)$ from $B\times F$ to the fiber product $B\times F$, which maps $(b,f)$ to $(b,g_{b}(f))$, where $g_{b}:F\cong \{b\}\times F\rightarrow \{g(b)\}\times F\cong F$, and the diagram
\[
\begin{CD}
B\times F @>\phi>> B\times F\\
@V\pi VV @VV\pi V\\
B @>id>> B
\end{CD}
\]
commutes. On the other hand, we have
\[
\begin{CD}
B\times F @>>> F\\
@V\pi VV @VVV\\
B @>>> \overline{\bbf}_{p}
\end{CD}
\]
Hence $R^{i}\pi_{!}(\mathbb{Z}/n\mathbb{Z})$ is the constant sheaf $H_{c}^{i}(F,\mathbb{Z}/n\mathbb{Z})$ on $B$. The automorphism $\phi$ acts on it and, at $b\in B$ it acts the way $g_{b}$ acts on $H_{c}^{i}(F,\mathbb{Z}/n\mathbb{Z})$. Since an endomorphism of a constant sheaf over a connected base is constant, the action of $\phi$ is the same everywhere. Passing to limit, we deduce that the actions of $g_{b}$ on $H_{c}^{*}(F,\mathbb{Q}_{l})$ are the same for every $b\in B$.

Suppose $E,B,F$ and the G-actions are defined over $\bbf_{q}$. Fix $n>0$. If $b_{1}, b_{2}\in B$ are fixed points of $gF_{q^{n}}$, then by what we just proved and the Lefschetz trace formula, we deduce that $|(\{b_1\}\times F)(\overline{\bbf}_{p})^{gF_{q^{n}}}|=|(\{b_{2}\}\times F)(\overline{\bbf}_{p})^{gF_{q^{n}}}|$. Hence we have the following equality.
\[
|(B\times F)(\overline{\bbf}_{p})^{gF_{q^{n}}}|=|B(\overline{\bbf}_{p})^{gF_{q^{n}}}||(\{b_1\}\times F)(\overline{\bbf}_{p})^{gF_{q^{n}}}|.
\]
Since the equality holds for all $n>0$, by the proof in Lemma \ref{Weil}, the lemma follows.

Now for the general case, we fix $g\in G$. It suffices to prove that the action of $g_{b_1}$ on $H_{c}^{*}(F,\mathbb{Q}_{l})$ is the same as the action of $g_{b_{2}}$ for any $b_{1}, b_{2}\in B$ fixed by $gF_{q^{n}}$. Take open neighborhoods $U_{1}, U_{2}$ of $b_{1}, b_{2}$ which trivialize the bundle. Replacing $U$ by $\cap_{n=0}^{\infty} g^{n}(U)$, we can assume $U_{1}, U_{2}$ are $g$-stable and connected since $B$ is irreducible. Now let $V=U_{1}\cap U_{2}$ and take any closed point $b_{0}\in V$. By the discussion in the trivial bundle case, we deduce that the action of $g_{b_1}$ is the same as the action of $g_{b_{0}}$, which is the same as the action of $g_{b_{2}}$. Hence we have 
\[
|E(\overline{\bbf}_{p})^{gF_{q^{n}}}|=|B(\overline{\bbf}_{p})^{gF_{q^{n}}}||(\{b_1\}\times F)(\overline{\bbf}_{p})^{gF_{q^{n}}}|
\]
for all $n>0$, and we are done.
\end{proof}

\begin{cor}\label{generalized jac}
Let $C$ be an integral projective curve over $\overline{\bbf}_{p}$ with an action of a finite group $G$. Then 
\[
[e(J^{n}C)]=[e(L)][e(J^{n}\tilde{C})]
\]
as virtual $G$-representations, i.e. in $R_{\mathbb{Q}_{l}}(G)$, where $L$ is a linear algebraic group and $\tilde{C}$ is the normalization of $C$.
\end{cor}

\begin{proof}
Let $f^{*}:J^{n}C\rightarrow J^{n}\tilde{C}$ be the pullback map. Since $g$ is an automorphism on $C$ and $\tilde{C}$, we have $g_{*}f^{*}=f^{*}g_{*}$. Now we use Lemma \ref{bundle}. 
\end{proof}

Now to prove Theorem \ref{curve1}, we first prove the following statement about $[e(J^{n}C)]$.

\begin{lemma}\label{trivial rep} 
Let $C$ be an integral curve over $\overline{\bbf}_{p}$ with nodal singularities and a $G$-action. Suppose $p\nmid|G|$. If $\tilde{C}/\left<g\right>$ is not a rational curve for every $g\in G$, then $[e({J^{n}C})]=0$ as $G$-representations.
\end{lemma}

\begin{proof}
By Corollary \ref{generalized jac}, it suffices to prove $[e(J^{n}\tilde{C})]=0$, which is equivalent to ${\Tr}(g,{[e(J^{n}\tilde{C})]})=0$ for any $g\in G$. But $J^{n}\tilde{C}$ is an abelian variety, which means $H^{i}(J^{n}\tilde{C},\mathbb{Q}_{l})\cong\wedge^{i}H^{1}(J^{n}\tilde{C},\mathbb{Q}_{l})$. Since $p\nmid|G|$, the curve $\tilde{C}/\left<g\right>$ is smooth. Then since $\tilde{C}/\left<g\right>$ is not rational, we have $H^{1}(\tilde{C}/\left<g\right>,\mathbb{Q}_{l})\neq 0$. Hence $H^{1}(J^{n}\tilde{C},\mathbb{Q}_{l})^{\left<g\right>}\neq 0$, which implies ${\Tr}(g,{[e(J^{n}\tilde{C})]})=0$. This is because $H^{1}(J^{n}\tilde{C},\mathbb{Q}_{l})=V_{0}\oplus V_{1}$, where $V_{0}$ is the non-empty eigenspace of $g$ with eigenvalue 1, and $V_{1}$ is its complement. It follows that $\sum(-1)^{i}[\wedge^{i}H^{1}(J^{n}\tilde{C},\mathbb{Q}_{l})]=(\sum(-1)^{i}[\wedge^{i}V_{0}])\otimes(\sum(-1)^{i}[\wedge^{i}V_{1}]$) and $\sum(-1)^{i}[\wedge^{i}V_{0}]$ has trace 0. 
\end{proof}

Now with the help of Corollary \ref{jac strata}, we can prove Theorem \ref{curve1} and Corollary \ref{curve2}.

\begin{proof}[Proof of Theorem \ref{curve1}]
Fix $g\in G$. Recall that $\overline{J^{n}C}=\coprod_{C'\rightarrow C} {J^{n}C'}$ by Corollary \ref{jac strata}. Depending on the action of $g$ on the nodes of $C$, $g_{*}$ permutes or acts on the strata $J^{n}C'$. For any union $\coprod J^{n}C'$ of two or more strata permuted by $g_{*}$ cyclically, the trace of $g$ on $H^{i}_{c}(\coprod J^{n}C')$ equals 0 since $g$ acts by cyclically permuting the components $H^{i}_{c}(J^{n}C')$. For the stratum which is stable under $g$, the trace of $g$ is also 0 by Lemma \ref{trivial rep}. Hence $[e(\overline{J^{n}C})]=0$ by Lemma \ref{strata}. 
\end{proof}

\begin{proof}[Proof of Corollary \ref{curve2}]
If $\tilde{C}/G$ is not a rational curve, then $\tilde{C}/\left<g\right>$ is not rational for any $g\in G$. 
\end{proof}

\section{Rational curves on surfaces}

Let $S$ be a smooth projective $K3$ surface over $\mathbb{C}$ with a $G$-action, and let $\mathcal{C}$ be the tautological family of curves over an $n$-dimensional integral linear system $|\mathcal{L}|$ acted on by $G$. Then $\overline{J^{n}\mathcal{C}}$ is a smooth projective variety over $|\mathcal{L}|$ whose fiber over a point $t\in\mathcal{L}$ is the compactified jacobian $\overline{J^{n}C_{t}}$. Choose some good reduction over $q$ such that `everything' ($\overline{J^{n}\mathcal{C}}$, $S$, $G$-action etc.) is defined over $\mathbb{F}_{q}$, where $q$ is a $p$ power for a prime number $p$. If we choose a large enough $p$, then $|\mathcal{L}|$ is still integral after the reduction. Notice that $|\mathcal{L}|$ has a stratification where each stratum $B$ satisfies ${\rm Stab}_{G}(t)=H$ for every $t\in B$ and some subgroup $H$, and the fibers $\mathcal{C}_{t}$ of the stratum have the same geometric genus. This is because for any subgroup $H$ in $G$, $|\mathcal{L}|^{H}\backslash\cup_{H'\supsetneqq H}|\mathcal{L}|^{H'}$ is a locally closed subspace. The reason for the stratification by the geometric genus is that the geometric genus gives a lower semicontinuous function in our case \cite[Proposition 2.4]{DH88}. Now notice that $gB=B$ if $g\in N_{G}(H)$ and $gB\cap B=\emptyset$ if $g\notin N_{G}(H)$. Hence we have a new stratification of $|\mathcal{L}|$ where each stratum $\cup_{g\in G}gB$ is $G$-invariant. 

\begin{proof}[Proof of Corollary \ref{euler}]
Recall that we have a birational map from $\overline{J^{n}\mathcal{C}}$ to $S^{[n]}$ \cite{Bea99}, which maps a pair $(C_{t},L)$ to a unique effective divisor $D$ on $C_{t}$ which is linearly equivalent to $L$, and can be viewed as a length $d$ subscheme of $S$. Note that the $G$-action commutes with this map. Hence by Corollary \ref{calabi-yau} $[e(\overline{J^{n}\mathcal{C}})]=[e(S^{[n]})]$ as $G$-representations. Now by Corollary \ref{trace} we are done. 
\end{proof}

Now we want to explicitly describe $\sum_{n=0}^{\infty}[e(S^{[n]})]t^{n}$ as a $G$-representation. For that purpose, we cite the following theorems (\cite[Proposition 1.2]{Mu88} and \cite[Lemma 2.3]{AST11}).

\begin{thm}\label{symplectic fact}
\cite{Mu88} Let $g$ be a symplectic automorphism of a complex K3 surface $S$ of order $n<\infty$. Then the number of fixed points of $g$ is equal to $\epsilon(n)=24\left(n\prod_{p|n}\left(1+\frac{1}{p}\right)\right)^{-1}$.
\end{thm}

\begin{thm}\label{nonsymplectic fact}
\cite{AST11} Let $g$ be a non-symplectic automorphism of a complex K3 surface $S$ of prime order $p$. Then the Euler characteristic of $S^{g}$ is $24-dp$, where $S^{g}$ denotes the fixed locus of $g$, $d=\frac{{\rm rank}\,T(g)}{p-1}$, and $T(g):=(H^{2}(S,\mathbb{Z})^{g})^{\bot}$.
\end{thm}

\begin{proof}[Proof of Theorem \ref{symplectic}]

By Corollary \ref{trace}, we deduce that 
\[
\sum_{n=0}^{\infty}\text{Tr}(g,[e({S}^{[n]})])t^{n}=\text{exp}\left(\sum_{m=1}^{\infty}\sum_{k=1}^{\infty}\frac{\text{Tr}(g^{k},[e({S})])t^{mk}}{k}\right).
\]
Then using the Lefschetz fixed point formula and Theorem \ref{symplectic fact}, we obtain the equality we want. 

When $G$ is a cyclic group of order $N$, we have $N\leq 8$ by \cite[Corollary 15.1.8]{H16}. Recall the definition of the Dedekind eta function $\eta(t)=t^{1/24}\prod_{n=1}^{\infty}(1-t^{n})$, where $t=e^{2\pi iz}$. Fix a generator $g$ of $G$.

If $N$ is a prime number $p$, we notice that ord$(g^{k})=1$ if $p|k$, and ord$(g^{k})=p$ otherwise. Hence
\[
\begin{aligned}
\sum_{n=0}^{\infty}\text{Tr}(g,[e({S}^{[n]})])t^{n}&=\text{exp}\left(\sum_{m=1}^{\infty}\sum_{k=1}^{\infty}\frac{\epsilon(p)t^{mk}}{k}\right)\text{exp}\left(\sum_{m=1}^{\infty}\sum_{k=1}^{\infty}\frac{(24-\epsilon(p))t^{mpk}}{pk}\right)\\
&=\left(\prod_{m=1}^{\infty}(1-t^{m})\right)^{-\epsilon(p)}\left(\prod_{m=1}^{\infty}(1-t^{mp})\right)^{\frac{\epsilon(p)-24}{p}}\\
&=\left(\prod_{m=1}^{\infty}(1-t^{m})(1-t^{mp})\right)^{-\frac{24}{p+1}}\\
&=t/\eta(t)^{\frac{24}{p+1}}\eta(t^{p})^{\frac{24}{p+1}}.
\end{aligned}
\]

If $N=4$, we have
\[
\text{Tr}(g^{k},[e({S})])=
\begin{cases}
4, & {\rm if}\ k\equiv 1,3\ ({\rm mod}\ 4)\\
8, & {\rm if}\ k\equiv 2\ ({\rm mod}\ 4)\\
24, & {\rm if}\ k\equiv 0\ ({\rm mod}\ 4).
\end{cases}
\]
Hence
\[
\begin{aligned}
\sum_{n=0}^{\infty}\text{Tr}(g,[e({S}^{[n]})])t^{n}&=\text{exp}\left(\sum_{m=1}^{\infty}\sum_{k\equiv 1,3}\frac{4t^{mk}}{k}+\sum_{m=1}^{\infty}\sum_{k\equiv 2}\frac{8t^{mk}}{k}+\sum_{m=1}^{\infty}\sum_{k\equiv 0}\frac{24t^{mk}}{k}\right)\\
&=\prod_{m=1}^{\infty}(1-t^{m})^{-4}(1-t^{2m})^{-2}(1-t^{4m})^{-4}\\
&=t/\eta(t)^{4}\eta(t^{2})^{2}\eta(t^{4})^{4}.
\end{aligned}
\]

If $N=6$, we have
\[
\text{Tr}(g^{k},[e({S})])=
\begin{cases}
2, & {\rm if}\ k\equiv 1,5\ ({\rm mod}\ 6)\\
6, & {\rm if}\ k\equiv 2,4\ ({\rm mod}\ 6)\\
8, & {\rm if}\ k\equiv 3\ ({\rm mod}\ 6)\\
24, & {\rm if}\ k\equiv 0\ ({\rm mod}\ 6).
\end{cases}
\]
Hence
\[
\begin{aligned}
\sum_{n=0}^{\infty}\text{Tr}(g,[e({S}^{[n]})])t^{n}&=\text{exp}\left(\sum_{m=1}^{\infty}\sum_{k\equiv 1,5}\frac{2t^{mk}}{k}+\sum_{m=1}^{\infty}\sum_{k\equiv 2,4}\frac{6t^{mk}}{k}+\sum_{m=1}^{\infty}\sum_{k\equiv 3}\frac{8t^{mk}}{k}+\sum_{m=1}^{\infty}\sum_{k\equiv 0}\frac{24t^{mk}}{k}\right)\\
&=\prod_{m=1}^{\infty}(1-t^{m})^{-2}(1-t^{2m})^{-2}(1-t^{3m})^{-2}(1-t^{6m})^{-2}\\
&=t/\eta(t)^{2}\eta(t^{2})^{2}\eta(t^{3})^{2}\eta(t^{6})^{2}.
\end{aligned}
\]

If $N=8$, we have
\[
\text{Tr}(g^{k},[e({S})])=
\begin{cases}
2, & {\rm if}\ k\equiv 1,3,5,7\ ({\rm mod}\ 8)\\
4, & {\rm if}\ k\equiv 2,6\ ({\rm mod}\ 8)\\
8, & {\rm if}\ k\equiv 4\ ({\rm mod}\ 8)\\
24, & {\rm if}\ k\equiv 0\ ({\rm mod}\ 8).
\end{cases}
\]
Hence
\[
\begin{aligned}
\sum_{n=0}^{\infty}\text{Tr}(g,[e({S}^{[n]})])t^{n}&=\text{exp}\left(\sum_{m=1}^{\infty}\sum_{k\equiv 1,3,5,7}\frac{2t^{mk}}{k}+\sum_{m=1}^{\infty}\sum_{k\equiv 2,6}\frac{4t^{mk}}{k}+\sum_{m=1}^{\infty}\sum_{k\equiv 4}\frac{8t^{mk}}{k}+\sum_{m=1}^{\infty}\sum_{k\equiv 0}\frac{24t^{mk}}{k}\right)\\
&=\prod_{m=1}^{\infty}(1-t^{m})^{-2}(1-t^{2m})^{-1}(1-t^{4m})^{-1}(1-t^{8m})^{-2}\\
&=t/\eta(t)^{2}\eta(t^{2})\eta(t^{4})\eta(t^{8})^{2}.
\end{aligned}
\]
\end{proof}

\begin{proof}[Proof of Theorem \ref{nonsymplectic}]

By Corollary \ref{trace}, we have
\[
\sum_{n=0}^{\infty}\text{Tr}(g,[e({S}^{[n]})])t^{n}=\text{exp}\left(\sum_{m=1}^{\infty}\sum_{k=1}^{\infty}\frac{\text{Tr}(g^{k},[e({S})])t^{mk}}{k}\right).
\]
By the topological Lefschetz formula, the Euler characteristic
\[
e(S^{g})=\sum_{i=0}^{4}(-1)^{i}{\rm Tr}(g^{*},H^{i}(S,\mathbb{C}))=\text{Tr}(g,[e({S})]).
\]
Fix $g\neq 1$ and notice that $S^{g}$ is the same as $S^{g^{k}}$ for $p\nmid k$. We deduce $\text{Tr}(g,[e({S})])=\text{Tr}(g^{k},[e({S})])=24-dp$ by Theorem \ref{nonsymplectic fact}. Hence
\[
\begin{aligned}
\sum_{n=0}^{\infty}\text{Tr}(g,[e({S}^{[n]})])t^{n}&=\text{exp}\left(\sum_{m=1}^{\infty}\sum_{k=1}^{\infty}\frac{(24-dp)t^{mk}}{k}\right)\text{exp}\left(\sum_{m=1}^{\infty}\sum_{k=1}^{\infty}\frac{(dp)t^{mpk}}{pk}\right)\\
&=\left(\prod_{m=1}^{\infty}(1-t^{m})\right)^{dp-24}\left(\prod_{m=1}^{\infty}(1-t^{mp})\right)^{-d}.
\end{aligned}
\]
\end{proof}

{\bf Example 1 ($\mathbb{Z}/2\mathbb{Z}$).} Here we look at an explicit K3 surface with a symplectic $\mathbb{Z}/2\mathbb{Z}$-action. Consider the elliptic K3 surface $S$ defined by the Weierstrass equation
\[
y^{2}=x^{3}+(t^{4}+a_{1}t^{2}+a_{2})x+(t^{12}+b_{1}t^{10}+b_{2}t^{8}+b_{3}t^{6}+b_{4}t^{4}+b_{5}t^{2}+b_{6}),
\]
where $(a_{1},a_{2})\in\mathbb{C}^{2}$, $(b_{1},...,b_{6})\in\mathbb{C}^{6}$ are generic. The fibration has 24 nodal fibers (Kodaira type $I_{1}$) over the zeros of its discriminant polynomial and those zeros do not contain $0$ and $\infty$. The automorphsim of order 2
\[
\sigma(x,y,t)=(x,-y,-t)
\]
acts non-trivially on the base of the fibration and preserves the smooth elliptic curves over $t=0$ and $t=\infty$. Now denote one of the fibers by $L$, then $|L|$ is a $\left<\sigma\right>$-invariant integral linear system and all of the singular curves in $|L|$ are nodal rational curves. We want to understand the $\sigma$-orbits in $|L|$.

Since we know explicitly the action of $\sigma$, by calculation we know that there are 4 $\sigma$-fixed points on the fiber over $t=0$ and 4 $\sigma$-fixed points on the fiber over $t=\infty$. So $\sigma$ has 8 isolated fixed points, hence it is a symplectic involution.

Now we know from Theorem \ref{symplectic} that
\[
\sum_{n=0}^{\infty}\text{Tr}(\sigma,[e({S}^{[n]})])t^{n}=\left(\prod_{k=1}^{\infty}(1-t^{k})\right)^{-8}\left(\prod_{k=1}^{\infty}(1-t^{2k})\right)^{-8}.
\] 
This implies ${\rm Tr}(\sigma,[e(S)])=8$ by looking at the coefficient of $t$. Hence we have ${\rm Tr}(\sigma,[e(\overline{J\mathcal{C}})])=8$ by Corollary \ref{euler}. But in this case, there are only two strata contributing to $[e(\overline{J\mathcal{C}})]$ as $\left<\sigma\right>$-representations (see Remark \ref{rmk}). One consists of nodal rational curves which are not $\sigma$-stable. The other consists of elliptic curves whose quotient by $\sigma$ is rational. The first stratum contributes $n_{1}[\mathbb{Z}/2\mathbb{Z}]$ if there are $n_{1}$ such $\sigma$-orbits in $|L|$, where $[\mathbb{Z}/2\mathbb{Z}]$ is the regular representation. The second stratum contributes $n_{2}(2V_{1}-2V_{-1})$ if there are $n_{2}$ such elliptic curves in $|L|$, where $V_{s}$ is the 1-dim representation on which $\sigma$ has eigenvalue $s$. This is because for a $\sigma$-stable smooth curve $C$ of genus g whose quotient by $\left<\sigma\right>$ is rational, by a similar argument as in the proof of Lemma \ref{trivial rep}, we deduce that $[e(J^{g}(C))]=2^{2g-1}V_{1}-2^{2g-1}V_{-1}$. Hence
\[
[e(\overline{J\mathcal{C}})]=n_{1}[\mathbb{Z}/2\mathbb{Z}]+n_{2}(2V_{1}-2V_{-1}).
\]
But since we already know the representation $[e(\overline{J\mathcal{C}})]$, by calculation we have $n_{1}=12$ and $n_{2}=2$.

On the other hand, this coincides with the geometric picture. From the definition of $\sigma$ we observe that there are indeed 12 $\sigma$-orbits of nodal rational curves. Denote by $C_{0}, C_{\infty}$ the fibers over $t=0, \infty$. Since $\sigma$ preserve $C_{0}$ and there are 4 $\sigma$-fixed points, we deduce that the degree 2 morphism $C_{0}\rightarrow C_{0}/\left<\sigma\right>$ has 4 ramification points. Hence by the Riemann-Hurwitz formula $C_{0}/\left<\sigma\right>$ is smooth rational. By the same argument, $C_{\infty}/\left<\sigma\right>$ is also smooth rational. This is what we have expected since there should be two such curves from the calculation of the representations.

{\bf Example 2 ($\mathbb{Z}/3\mathbb{Z}$).} Here we look at an explicit K3 surface with a non-symplectic $\mathbb{Z}/3\mathbb{Z}$-action \cite[Remark 4.2]{AST11}. Consider the elliptic K3 surface $S$ defined by the Weierstrass equation
\[
y^{2}=x^{3}+(t^6+a_{1}t^{3}+a_{2})x+(t^{12}+b_{1}t^{9}+b_{2}t^{6}+b_{3}t^{3}+b_{4}),
\]
where $(a_{1},a_{2})\in\mathbb{C}^{2}$, $(b_{1},...,b_{4})\in\mathbb{C}^{4}$ are generic. The fibration has 24 nodal fibers (Kodaira type $I_{1}$) over the zeros of its discriminant polynomial and those zeros do not contain $0$ and $\infty$. The automorphism of order 3
\[
\sigma(x,y,t)=(x,y,\zeta_{3}t)
\]
acts non-trivially on the basis of the fibration and preserves the smooth elliptic curves over $t=0$ and $t=\infty$. Now denote one of the fibers by $L$, then $|L|$ is a $\left<\sigma\right>$-invariant integral linear system and all of the singular curves in $|L|$ are nodal rational curves. We want to understand the $\sigma$-orbits in $|L|$.

First we observe that $\sigma$ fixs the fiber over $t=0$. Hence $\left<\sigma\right>$ acts non-symplectically on $S$. We know from \cite[Theorem 4.1]{AST11} that rank $T(\sigma)=14$. Then using Theorem \ref{nonsymplectic} we have 
\[
\sum_{n=0}^{\infty}\text{Tr}(\sigma,[e({S}^{[n]})])t^{n}=\left(\prod_{k=1}^{\infty}(1-t^{k})\right)^{-3}\left(\prod_{k=1}^{\infty}(1-t^{3k})\right)^{-7}.
\]
This implies ${\rm Tr}(\sigma,[e(S)])=3$ by looking at the coefficient of $t$. Hence we have ${\rm Tr}(\sigma,[e(\overline{J\mathcal{C}})])=3$ by Corollary \ref{euler}. But in this case, there are only two strata contributing to $[e(\overline{J\mathcal{C}})]$ as $\left<\sigma\right>$-representations (see Remark \ref{rmk}). One consists of nodal rational curves which are not $\sigma$-stable. The other consists of elliptic curves whose quotient by $\sigma$ is rational. In particular, the fiber over $t=0$ will not contribute to $[e(\overline{J\mathcal{C}})]$. The first stratum contributes $n_{1}[\mathbb{Z}/3\mathbb{Z}]$ if there are $n_{1}$ such $\sigma$-orbits in $|L|$, where $[\mathbb{Z}/3\mathbb{Z}]$ is the regular representation. The second stratum contributes $n_{2}(2V_{1}-V_{\zeta_{3}}-V_{\zeta^{-1}_{3}})$ if there are $n_{2}$ such elliptic curves in $|L|$ , where $V_{s}$ is the 1-dim representation on which $\sigma$ has eigenvalue $s$. Hence
\[
[e(\overline{J\mathcal{C}})]=n_{1}[\mathbb{Z}/3\mathbb{Z}]+n_{2}(2V_{1}-V_{\zeta_{3}}-V_{\zeta^{-1}_{3}}).
\]
But since we already know the representation $[e(\overline{J\mathcal{C}})]$, by calculation we have $n_{1}=8$ and $n_{2}=1$.

On the other hand, this coincides with the geometric picture. From the definition of $\sigma$ we observe that there are indeed 8 $\sigma$-orbits of nodal rational curves. Since the action of $\sigma$ is explicit, by calculation we know that there are 3 $\sigma$-fixed points on the fiber over $t=\infty$. Denote by $C_{\infty}$ the fiber over $t=\infty$. Then this implies that the degree 3 morphism $C_{\infty}\rightarrow C_{\infty}/\left<\sigma\right>$ has 3 ramification points each of order 2. Hence by the Riemann-Hurwitz formula $C_{\infty}/\left<\sigma\right>$ is smooth rational, which is what we have expected.

\begin{remark}{\bf ($\mathbb{Z}/2\mathbb{Z}$ in general)}\label{rmk}
Let us consider a complex K3 surface $S$ with a $\mathbb{Z}/2\mathbb{Z}$-action (i.e. an involution $\sigma$). Take a $\sigma$-invariant integral linear system $|\mathcal{L}|$ of dimension $d$. Assume all the rational curves in $|\mathcal{L}|$ have nodal singularities and $|\mathcal{L}|$ has finitely many $\sigma$-fixed points.

For the stratum of $|\mathcal{L}|$ which consists of curves that are not $\sigma$-stable and of geometric genus $>0$, denote by $M$ the corresponding stratum of $\overline{J^{d}\mathcal{C}}$. Then $[e(M)]=e(M/\left<\sigma\right>)[\mathbb{Z}/2\mathbb{Z}]=e(\overline{J^{d}C_{0}})e(B/\left<\sigma\right>)[\mathbb{Z}/2\mathbb{Z}]$, where $[\mathbb{Z}/2\mathbb{Z}]$ is the regular representation. But $e(\overline{J^{d}C_{0}}):=\sum_{n=0}^{\infty} (-1)^{n}{\rm dim}H^{n}(\overline{J^{d}C_{0}},\mathbb{Q}_{l})=0$ by Theorem \ref{curve1}. Hence $[e(M)]=0$ as representation.

For the stratum of $|\mathcal{L}|$ which consists of curves that are not $\sigma$-stable and are nodal rational curves, we deduce $[e(M)]=n_{0}[\mathbb{Z}/2\mathbb{Z}]$, where $n_{0}$ is the number of $\sigma$-orbits of nodal rational curves, since $e(\overline{J^{d}C_{0}})=1$ by Corollary \ref{jac strata}, Corollary \ref{generalized jac} and $e(\mathbb{G}_{m})=0$.

For those curves that are $\sigma$-stable, we first notice that $\sigma$ fixes some curve only if $\sigma$ acts non-symplectically on $S$, and in that case the fixed curves are always smooth (\cite[$\S 2$]{AST11}). Let $C_{0}$ be a smooth curve of genus $d\geq 1$ fixed by $\sigma$. If $d=1$, then the fixed locus of $\sigma$ consists of two disjoint elliptic curves and they are linearly equivalent. If $d>1$, then it is the only curve of genus $d$ fixed by $\sigma$ (\cite[Theorem 3.1]{AST11}). In either case, the stratum consisting of $\sigma$-fixed curves contributes 0 to the representation since the Euler characteristic of an abelian variety is $0$.

For the stratum of $|\mathcal{L}|$ which consists of curves that are $\sigma$-stable, if the normalization of the curves quotient by $\left<\sigma\right>$ is not rational, then by Theorem 1.9, we have $[e(M)]=0$.

Hence there are only two strata contributing to $[e(\overline{J^{d}\mathcal{C}})]$. One consists of $\sigma$-orbits of nodal rational curves. The other consists of the curves whose normalization quotient by $\left<\sigma\right>$ are smooth rational.
\end{remark}
Now let us give an example of a 2-dim $\mathbb{Z}/2\mathbb{Z}$-invariant linear system. This example is suggested by Jim Bryan.

{\bf Example 3 (2-dim $\mathbb{Z}/2\mathbb{Z}$).} Let $S$ be a K3 surface given by the double cover of $\mathbb{P}^{2}$ branched over a smooth sextic curve $C$ in $\mathbb{P}^{2}$. Let $\tau$ be the involution on $\mathbb{P}^{2}$ sending $(x: y: z)$ to $(-x: y: z)$. Denote the covering involution by $i: S\rightarrow S$. Then if we suppose $C$ is $\tau$-invariant, the `composition' of $\tau$ and $i$ will give a symplectic involution $\sigma$ on $S$ (\cite[Section 3.2]{GS07}). 

The fixed locus of $\tau$ on $\mathbb{P}^{2}$ consists of a point $x_{0}=(1: 0: 0)$ and a line $l_{0}=\{(x: y: z)| x=0\}$. Denote the six intersection points of $l_{0}$ and $C$ by $x_{3}, x_{4},..., x_{8}$. Let $\pi: S\rightarrow\mathbb{P}^{2}$ be the double cover map. Denote the two points in $\pi^{-1}(x_{0})$ by $x_{1}, x_{2}$. Then the fixed locus of $\sigma$ is the eight points $x_{1}, x_{2},..., x_{8}$. Notice that $\sigma$ commutes with $i$ and the induced action of $\sigma$ on $\mathbb{P}^{2}$ is just $\tau$.

Now let $\mathcal{L}=\pi^{*}\mathcal{O}_{\mathbb{P}^{2}}(1)$. Then the linear system $|\mathcal{L}|$ consists of the curves which are the preimages of the lines in $\mathbb{P}^{2}$ under $\pi$. For a generic choice of $C$, $|\mathcal{L}|$ is a $\sigma$-invariant integral linear system. A generic line will intersects $C$ in six points, and its preimage is a smooth genus 2 curve. Some lines will intersect $C$ in a tangent point and 4 other distinct points, and their preimages are curves with one node. The other lines are the 324 bitangents of $C$, which can be seen from the Pl\"ucker formula or the coefficient of $t^{2}$ in $\prod_{n=1}^{\infty}(1-t^{n})^{-24}$.

Let $\mathcal{C}\rightarrow|\mathcal{L}|$ be the tautological family of curves over $|\mathcal{L}|$. Now we know from Theorem \ref{symplectic} that
\[
\sum_{n=0}^{\infty}\text{Tr}(\sigma,[e({S}^{[n]})])t^{n}=\left(\prod_{k=1}^{\infty}(1-t^{k})\right)^{-8}\left(\prod_{k=1}^{\infty}(1-t^{2k})\right)^{-8}.
\] 

This implies ${\rm Tr}(\sigma,[e({S}^{[2]})])=52$ by looking at the coefficient of $t^{2}$. Hence we have ${\rm Tr}(\sigma,[e(\overline{J\mathcal{C}})])=52$ by Corollary \ref{euler}. Since we know ${\rm Tr}(1,[e(S^{[2]})])=324$, we obtain
\[
[e(\overline{J\mathcal{C}})]=188V_{1}+136V_{-1},
\]
where $V_{1}$ is the 1-dimensional trivial representation and $V_{-1}$ is the 1-dimensional representation on which $\sigma$ has eigenvalue $-1$.

To give a geometric interpretation of this representation, we first concentrate on the $\sigma$-invariant curves in $|\mathcal{L}|$ since non-invariant curves with zero or one node will contribute nothing to the representation, and non-invariant nodal rational curves will contribute some multiple of the regular representation. The $\sigma$-invariant curves in $|\mathcal{L}|$ are given by the preimages of the $\tau$-invariant lines, which consists of the line $l_{0}=\{(x: y: z)| x=0\}$ and all the lines passing through $x_{0}=(1: 0: 0)$, i.e., $\{(x: y: z)| by+cz=0\}$, $(b: c)\in\mathbb{P}^{1}$. 

The preimage of $l_{0}$ is a smooth genus 2 curve, and it has 6 ramification points $x_{3}, x_{4},..., x_{8}$ under $\sigma$. Hence its preimage quotient by $\sigma$ is a smooth rational curve, and it will contribute $8V_{1}-8V_{-1}$ to the representation by a similar argument as in the proof of Lemma \ref{trivial rep}.

The preimages of $\{(x: y: z)| by+cz=0\}$ are more complicated. A generic line will intersect $C$ in six points, and its preimage is a smooth curve of genus 2. It has 2 ramification points $x_{1}, x_{2}$ under $\sigma$. Hence its preimage quotient by $\sigma$ is an elliptic curve, and it will contribute nothing to the representation. 

Now let us consider those tangent lines. We first observe that if a line passes through one of the six points $x_{3}, x_{4},..., x_{8}$, then by looking at the $\sigma$-action on the preimage of the line, this point must be a tangent point. So its preimage is a curve with one node, and the normalization of it has 4 ramification points $x_{1}, x_{2}, x_{i}^{1}, x_{i}^{2}$ under $\sigma$, where $x_{i}^{1}, x_{i}^{2}$ are the two points on the normalization over the point $x_{i}$ if our line passes through $x_{i}$, $i\in\{3,4,...,8\}$. Hence the normalization of its preimage quotient by $\sigma$ is a rational curve, and we denote its contribution to the representation by $[e(\overline{J{C}_{1}})]$.

Now if there is a tangent point $y$ of the tangent line which is not one of the six points $x_{3}, x_{4},..., x_{8}$, then since both of the line and the curve $C$ are $\tau$-invariant, this line must have another tangent point $\tau(y)$. Hence this line must be a bitangent. In order to calculate the number of the $\tau$-invariant bitangents, we notice that the degree of the dual curve of $C$ is 30. We already have 6 tangent lines with one tangent point, and all the other $\tau$-invariant tangent lines are bitangents. Hence there should be 12 $\tau$-invariant bitangents. The preimage of the $\tau$-invariant bitangent is a rational curve with two nodes, and $\tau$ permutes the nodes. The normalization of the preimage has 2 ramification points $x_{1}, x_{2}$ under $\sigma$. Hence the normalization of its preimage quotient by $\sigma$ is a rational curve, and we denote its contribution to the representation by $[e(\overline{J{C}_{2}})]$. 

Finally, non-invariant curves with two nodes will also contribute to the representation. We know there are 324 curves with two nodes, and 12 of them are $\sigma$-invariant by the discussion above. So there are 312 non-invariant nodal rational curves, which will contribute $156[\mathbb{Z}/2\mathbb{Z}]=156V_{1}+156V_{-1}$ to the representation.

Combining all of the above, we deduce
\[
\begin{aligned}
188V_{1}+136V_{-1}&=[e(\overline{J\mathcal{C}})]\\
&=8V_{1}-8V_{-1}+6[e(\overline{J{C}_{1}})]+12[e(\overline{J{C}_{2}})]+156V_{1}+156V_{-1}.
\end{aligned}
\]
Hence we only need to check
\[
[e(\overline{J{C}_{1}})]+2[e(\overline{J{C}_{2}})]=4V_{1}-2V_{-1}
\]
and this is true by the following two lemmas.

\begin{lemma}\label{lemma1}
Let $C_{1}$ be an integral curve of arithmetic genus 2 with one node over $\overline{\bbf}_{p}$. If there is a involution $\sigma$ acting on it, and the action $\sigma$ on its normalization $\tilde{C_{1}}$ has 4 fixed points, two of which are the points over the node, then we have 
\[
[e(\overline{J{C}_{1}})]=2V_{1}-2V_{-1}
\]
as $\mathbb{Z}/2\mathbb{Z}$-representations.
\end{lemma}

\begin{lemma}\label{lemma2}
Let $C_{2}$ be an integral curve of arithmetic genus 2 with two nodes over $\overline{\bbf}_{p}$. If there is a involution $\sigma$ acting on it, and $\sigma$ permutes the nodes, then we have 
\[
[e(\overline{J{C}_{2}})]=V_{1}
\]
as $\mathbb{Z}/2\mathbb{Z}$-representations.
\end{lemma}

\begin{proof}[Proof of Lemma \ref{lemma1}]

Since $C_{1}$ is an integral curve of arithmetic genus 2 with one node, its normalization $\pi:\tilde{C_{1}}\rightarrow C_{1}$ is an elliptic curve, and we denote it by $E$. By Corollary 5.3 and Corollary 5.6, we have
\[
[e(\overline{J{C}_{1}})]=[e({J{C}_{1}})]+[e({J{E}})]=[e(\mathbb{G}_{m})][e(JE)]+[e(JE)].
\]

Now we notice that $\sigma$ has 4 fixed points on $E$, so $E\rightarrow E/\left<\sigma\right>$ realizes $E$ as a double cover of $\mathbb{P}^{1}$. Hence 
\[
[e(JE)]=[H^{0}(E,\mathbb{Q}_{l})]-[H^{1}(E,\mathbb{Q}_{l})]+[H^{2}(E,\mathbb{Q}_{l})]=2V_{1}-2V_{-1}.
\]

On the other hand, since $\sigma$ fixes the points over the node, we deduce from the short exact sequence 
\[
0\rightarrow\mathcal{O}_{C_{1}}^{*}\rightarrow\pi_{*}\mathcal{O}_{E}^{*}\rightarrow\delta\rightarrow 0
\]
that $\sigma$ acts trivially on the skyscraper sheaf $\delta$. Hence $\sigma$ acts trivially on $H^{0}(C_{1}, \delta)=\mathbb{G}_{m}$. So $[e(\mathbb{G}_{m})]=e(\mathbb{G}_{m})V_{1}=0$ since the topological Euler characteristic of $\mathbb{G}_{m}$ is 0.

Combining the above discussion, we have
\[
[e(\overline{J{C}_{1}})]=2V_{1}-2V_{-1}.
\]
\end{proof}

\begin{proof}[Proof of Lemma \ref{lemma2}]

Since $C_{2}$ is an integral curve of arithmetic genus 2 with two nodes, its normalization $\pi:\tilde{C_{2}}\rightarrow C_{2}$ is a rational curve. It also has two partial normalizations by resolving one of the nodes $\pi_{1}:C_{2}^{'}\rightarrow C_{2}$ and $\pi_{2}:C_{2}^{''}\rightarrow C_{2}$. By Corollary 5.3 and Corollary 5.6, we have
\[
[e(\overline{J{C}_{2}})]=[e({J{C}_{2}})]+[e({J{C}^{'}_{2}})]+[e({J{C}^{''}_{2}})]+[e({J\mathbb{P}^{1}})]=[e(\mathbb{G}_{m}\times\mathbb{G}_{m})]+e(\mathbb{G}_{m})[\mathbb{Z}/2\mathbb{Z}]+V_{1}
\]
since $\sigma$ permutes two nodes.

On the other hand, $\mathbb{G}_{m}$ is an affine curve. So dim $H^{2}_{c}(\mathbb{G}_{m},\mathbb{Q}_{l})=1$ and dim $H^{0}_{c}(\mathbb{G}_{m},\mathbb{Q}_{l})=0$. Since the topological Euler characteristic of $\mathbb{G}_{m}$ is 0, we also have dim $H^{1}_{c}(\mathbb{G}_{m},\mathbb{Q}_{l})=1$. Notice that $\sigma$ permutes two $\mathbb{G}_{m}$'s, and hence by the K\"unneth formula we have
\[
\begin{aligned}{}
[e(\mathbb{G}_{m}\times\mathbb{G}_{m})]&=[H^{1}_{c}(\mathbb{G}_{m},\mathbb{Q}_{l})\otimes H^{1}_{c}(\mathbb{G}_{m},\mathbb{Q}_{l})]-[H^{1}_{c}(\mathbb{G}_{m},\mathbb{Q}_{l})\otimes H^{2}_{c}(\mathbb{G}_{m},\mathbb{Q}_{l})]\\
&-[H^{2}_{c}(\mathbb{G}_{m},\mathbb{Q}_{l})\otimes H^{1}_{c}(\mathbb{G}_{m},\mathbb{Q}_{l})]+[H^{2}_{c}(\mathbb{G}_{m},\mathbb{Q}_{l})\otimes H^{2}_{c}(\mathbb{G}_{m},\mathbb{Q}_{l})]\\
&=V_{-1}-[\mathbb{Z}/2\mathbb{Z}]+V_{1}\\
&=0.
\end{aligned}
\]

Combining the above discussion, we have
\[
[e(\overline{J{C}_{2}})]=V_{1}.
\]
\end{proof}

Finally, let us give some discussions when $G$ equals a certain finite simple group.

{\bf Example 4 (PSL(2,7)).} Let $S$ be a complex K3 surface acting faithfully by $G=PSL_{2}(\mathbb{F}_{7})$. Such a K3 surface exists. For example, $PSL(2,7)$ acts faithfully and symplectically on the surface $X^{3}Y+Y^{3}Z+Z^{3}X+T^{4}=0$ in $\mathbb{P}^{3}$ by means of a linear action on $\mathbb{P}^{3}$ \cite{Mu88}. We know from Theorem \ref{curve1} that a $G$-stable curve $C$ with nodal singularities in an integral linear system does contribute to the representation $[e(\overline{J^{d}\mathcal{C}})]$ only if there exists some $g\in G$ such that $\tilde{C}/\left<g\right>=\mathbb{P}^{1}$. It turns out that if this happens, then $\tilde{C}$ must be the Klein quartic, which is the Hurwitz surface of the lowest possible genus. Notice that $G$ acts on $C$ faithfully since any non-trivial element of $G$ acts symplectically on $S$ and cannot fix curves.

\begin{prop}
Let $C$ be a smooth projective curve over $\mathbb{C}$ with a faithful $G=PSL(2,7)$-action. If there exists $g\in PSL(2,7)$ such that ${C}/\left<g\right>=\mathbb{P}^{1}$, then the genus of $C$ is $3$ and $g$ has order $7$. In particular, the automorphism group of $C$ reaches its Hurwitz bound, and hence $C$ is the Klein quartic.
\end{prop}

\begin{proof}
The idea is to use the equivariant Riemann-Hurwitz formula \cite[Chapter VI $\S 4$]{Se79} for $\pi:C\rightarrow C/G=\mathbb{P}^{1}$. We have 
\[
[e(C)]=e(\mathbb{P}^{1})I_{\left<1\right>}-\sum_{p\in\mathbb{P}^{1}}(I_{\left<1\right>}-I_{\left<h_{p}\right>})
\]
as $G$-representations, where $\left<h_{p}\right>$ is the stablizer of some point over $p$, $I_{\left<h_{p}\right>}$ denotes the induced representation ${\rm Ind}^{G}_{\left<h_{p}\right>}\mathbbm{1}$, and $\mathbbm{1}$ is the 1-dim trivial representation. Notice that $I_{\left<h_{p}\right>}$ is independent of the point we choose over $p$.

Since $G$ acts trivially on $\mathbb{P}^{1}$, we have
\[
H^{1}(C,\mathbb{C})=\sum_{p\in\mathbb{P}^{1}}(I_{\left<1\right>}-I_{\left<h_{p}\right>})-2I_{\left<1\right>}+2\mathbbm{1}
\]
and what we are going to do is to compare the representations on both sides. For this purpose, we need the character table of $PSL(2,7)$.
\begin{center}
\begin{tabular}{c||cccccc}
& $1A_{1}$ & $2A_{21}$ & $3A_{56}$ & $4A_{42}$ & $7A_{24}$ & $7B_{24}$\\
\hline\hline
$\chi_{1}$ & 1&1&1&1&1&1\\
$\chi_{2}$ & $1^{(3)}$&$1(-1)^{(2)}$&$1\omega\bar{\omega}$&1i(-i)&$\zeta_{7}\zeta_{7}^{2}\zeta_{7}^{4}$&$\zeta_{7}^{3}\zeta_{7}^{5}\zeta_{7}^{6}$\\
$\chi_{3}$&$1^{(3)}$& $1(-1)^{(2)}$&$1\omega\bar{\omega}$&1i(-i)&$\zeta_{7}^{3}\zeta_{7}^{5}\zeta_{7}^{6}$&$\zeta_{7}\zeta_{7}^{2}\zeta_{7}^{4}$\\
$\chi_{4}$&$1^{(6)}$& $1^{(4)}(-1)^{(2)}$&$1^{(2)}\omega^{(2)}\bar{\omega}^{(2)}$&$1^{(2)}(-1)^{(2)}i(-i)$&$\zeta_{7}...\zeta_{7}^{6}$&$\zeta_{7}...\zeta_{7}^{6}$\\
$\chi_{5}$&$1^{(7)}$& $1^{(3)}(-1)^{(4)}$&$1^{(3)}\omega^{(2)}\bar{\omega}^{(2)}$&$1(-1)^{(2)}i^{(2)}(-i)^{(2)}$&$1\zeta_{7}...\zeta_{7}^{6}$&$1\zeta_{7}...\zeta_{7}^{6}$\\
$\chi_{6}$&$1^{(8)}$& $1^{(4)}(-1)^{(4)}$&$1^{(2)}\omega^{(3)}\bar{\omega}^{(3)}$&$1^{(2)}(-1)^{(2)}i^{(2)}(-i)^{(2)}$&$1^{(2)}\zeta_{7}...\zeta_{7}^{6}$&$1^{(2)}\zeta_{7}...\zeta_{7}^{6}$
\end{tabular}
\end{center}

This is a refined character table which can be deduced from the usual character table. Each entry denotes the eigenvalues of the element in given conjuacy classes acting on given irreducible representations, $nA_{m}$ denotes the conjugacy class of size $m$ in which each element has order $n$, $a^{(i)}b^{(j)}$ denotes eigenvalue $a$ with multiplicity $i$ and eigenvalue $b$ with multiplicity $j$, and $\zeta_{7}...\zeta_{7}^{6}$ means $\zeta_{7}\zeta_{7}^{2}\zeta_{7}^{3}\zeta_{7}^{4}\zeta_{7}^{5}\zeta_{7}^{6}$.

For induced representations, the character $I_{\left<h_{p}\right>}(x)=\frac{1}{|\left<h_{p}\right>|}\sum_{u\in G}\chi(uxu^{-1})$, where $\chi(x)=1$ if $x\in\left<h_{p}\right>$, and $\chi(x)=0$ otherwise. Hence using our character table and calculating by Schur orthogonality relations, we have
\[
\begin{aligned}
I_{2}&=\chi_{1}+\chi_{2}+\chi_{3}+4\chi_{4}+3\chi_{5}+4\chi_{6}\\
I_{3}&=\chi_{1}+\chi_{2}+\chi_{3}+2\chi_{4}+3\chi_{5}+2\chi_{6}\\
I_{4}&=\chi_{1}+\chi_{2}+\chi_{3}+2\chi_{4}+\chi_{5}+2\chi_{6}\\
I_{7}&=\chi_{1}+\chi_{5}+2\chi_{6}
\end{aligned}
\]
where $I_{n}$ denotes $I_{\left<h_{p}\right>}$ for the element $h_{p}$ of order $n$.

Now since there exists $g\in G$ such that ${C}/\left<g\right>=\mathbb{P}^{1}$, we have $H^{1}(C,\mathbb{C})^{g}=H^{1}(C/\left<g\right>,\mathbb{C})=0$. But for $\chi_{1}$, $\chi_{5}$ and $\chi_{6}$, whatever $g$ is, there are always non-trivial $g$-fixed vectors. This implies that $H^{1}(C,\mathbb{C})$ does not contain $\chi_{1}$, $\chi_{5}$ and $\chi_{6}$ at all. Hence the coefficients of $\chi_{1}$, $\chi_{5}$ and $\chi_{6}$ in $\sum_{p\in\mathbb{P}^{1}}(I_{\left<1\right>}-I_{\left<h_{p}\right>})-2I_{\left<1\right>}+2\mathbbm{1}$ must be 0. This gives us only two possibilities: $H^{1}(C,\mathbb{C})=I_{1}-I_{2}-I_{3}-I_{4}+2\mathbbm{1}$ or $H^{1}(C,\mathbb{C})=I_{1}-I_{2}-I_{3}-I_{7}+2\mathbbm{1}$. If we look at the dimensions of the right hand sides, the first one gives dimension -12 and the second gives dimension 6. Hence the only possibility is $H^{1}(C,\mathbb{C})=I_{1}-I_{2}-I_{3}-I_{7}+2\mathbbm{1}=\chi_{2}+\chi_{3}$, which shows that the genus of $C$ is $\frac{1}{2}\text{dim}H^{1}(C,\mathbb{C})=3$. We also deduce from this argument that $g$ must has order 7 since the element of order not equal to 7 does have fixed vectors in $\chi_{2}$ and $\chi_{3}$.
\end{proof}

Following this observation, we do the calculations for some other groups in Mukai's list \cite{Mu88}.

{\bf Example 5 ($A_6$).} $G=A_{6}$ acts faithfully and symplectically on the K3 surface $\sum_{1}^{6}X_{i}=\sum_{1}^{6}X^{2}_{i}=\sum_{1}^{6}X^{3}_{i}=0$ in $\mathbb{P}^{5}$ via permutation action of coordinates on $\mathbb{P}^{5}$. Then by Theorem \ref{curve1}, a $G$-stable curve $C$ with nodal singularities in an integral linear system will not contribute to the representation $[e(\overline{J^{d}\mathcal{C}})]$.

\begin{prop}
Let $C$ be a smooth projective curve over $\mathbb{C}$ with a faithful $G=A_{6}$-action. Then for any $g\in A_{6}$, we have ${C}/\left<g\right>\neq\mathbb{P}^{1}$.
\end{prop}

\begin{proof}
We have the following character table for $A_{6}$.

\begin{center}
\begin{tabular}{c||cccc}
& $1A_{1}$ & $2A_{45}$ & $3A_{40}$ & $3B_{40}$\\
\hline\hline
$\chi_{1}$ & 1&1&1&1\\
$\chi_{2}$ & $1^{(5)}$&$1^{(3)}(-1)^{(2)}$&$1^{(3)}\omega\bar{\omega}$&$1\omega^{(2)}\bar{\omega}^{(2)}$\\
$\chi_{3}$&$1^{(5)}$& $1^{(3)}(-1)^{(2)}$&$1\omega^{(2)}\bar{\omega}^{(2)}$&$1^{(3)}\omega\bar{\omega}$\\
$\chi_{4}$&$1^{(8)}$& $1^{(4)}(-1)^{(4)}$&$1^{(2)}\omega^{(3)}\bar{\omega}^{(3)}$&$1^{(2)}\omega^{(3)}\bar{\omega}^{(3)}$\\
$\chi_{5}$&$1^{(8)}$& $1^{(4)}(-1)^{(4)}$&$1^{(2)}\omega^{(3)}\bar{\omega}^{(3)}$&$1^{(2)}\omega^{(3)}\bar{\omega}^{(3)}$\\
$\chi_{6}$&$1^{(9)}$& $1^{(5)}(-1)^{(4)}$&$1^{(3)}\omega^{(3)}\bar{\omega}^{(3)}$&$1^{(3)}\omega^{(3)}\bar{\omega}^{(3)}$\\
$\chi_{7}$&$1^{(10)}$& $1^{(4)}(-1)^{(6)}$&$1^{(4)}\omega^{(3)}\bar{\omega}^{(3)}$&$1^{(4)}\omega^{(3)}\bar{\omega}^{(3)}$
\end{tabular}
\end{center}
\begin{center}
\begin{tabular}{c||ccc}
& $4A_{90}$ & $5A_{72}$ & $5B_{72}$\\
\hline\hline
$\chi_{1}$ & 1&1&1\\
$\chi_{2}$ & $1{(-1)}^{(2)}i(-i)$&$1\zeta_{5}\zeta_{5}^{2}\zeta_{5}^{3}\zeta_{5}^{4}$&$1\zeta_{5}\zeta_{5}^{2}\zeta_{5}^{3}\zeta_{5}^{4}$\\
$\chi_{3}$ & $1{(-1)}^{(2)}i(-i)$&$1\zeta_{5}\zeta_{5}^{2}\zeta_{5}^{3}\zeta_{5}^{4}$&$1\zeta_{5}\zeta_{5}^{2}\zeta_{5}^{3}\zeta_{5}^{4}$\\
$\chi_{4}$ & $1^{(2)}(-1)^{(2)}i^{(2)}(-i)^{(2)}$&$1^{(2)}\zeta_{5}\zeta_{5}^{2(2)}\zeta_{5}^{3(2)}\zeta_{5}^{4}$&$1^{(2)}\zeta_{5}^{(2)}\zeta_{5}^{2}\zeta_{5}^{3}\zeta_{5}^{4(2)}$\\
$\chi_{5}$ & $1^{(2)}(-1)^{(2)}i^{(2)}(-i)^{(2)}$&$1^{(2)}\zeta_{5}^{(2)}\zeta_{5}^{2}\zeta_{5}^{3}\zeta_{5}^{4(2)}$&$1^{(2)}\zeta_{5}\zeta_{5}^{2(2)}\zeta_{5}^{3(2)}\zeta_{5}^{4}$\\
$\chi_{6}$ & $1^{(3)}(-1)^{(2)}i^{(2)}(-i)^{(2)}$&$1\zeta_{5}^{(2)}\zeta_{5}^{2(2)}\zeta_{5}^{3(2)}\zeta_{5}^{4(2)}$&$1\zeta_{5}^{(2)}\zeta_{5}^{2(2)}\zeta_{5}^{3(2)}\zeta_{5}^{4(2)}$\\
$\chi_{7}$ & $1^{(2)}(-1)^{(2)}i^{(3)}(-i)^{(3)}$&$1^{(2)}\zeta_{5}^{(2)}\zeta_{5}^{2(2)}\zeta_{5}^{3(2)}\zeta_{5}^{4(2)}$&$1^{(2)}\zeta_{5}^{(2)}\zeta_{5}^{2(2)}\zeta_{5}^{3(2)}\zeta_{5}^{4(2)}$
\end{tabular}
\end{center}

For induced representations, we have
\[
\begin{aligned}
I_{2}&=\chi_{1}+3\chi_{2}+3\chi_{3}+4\chi_{4}+4\chi_{5}+5\chi_{6}+4\chi_{7}\\
I_{3A}&=\chi_{1}+3\chi_{2}+\chi_{3}+2\chi_{4}+2\chi_{5}+3\chi_{6}+4\chi_{7}\\
I_{3B}&=\chi_{1}+\chi_{2}+3\chi_{3}+2\chi_{4}+2\chi_{5}+3\chi_{6}+4\chi_{7}\\
I_{4}&=\chi_{1}+\chi_{2}+\chi_{3}+2\chi_{4}+2\chi_{5}+3\chi_{6}+2\chi_{7}\\
I_{5}&=\chi_{1}+\chi_{2}+\chi_{3}+2\chi_{4}+2\chi_{5}+\chi_{6}+2\chi_{7}
\end{aligned}
\]

Now suppose there exists $g\in G$ such that ${C}/\left<g\right>=\mathbb{P}^{1}$. Then we have $H^{1}(C,\mathbb{C})^{g}=H^{1}(C/\left<g\right>,\mathbb{C})=0$. But for all the irreducible representations of $G$, whatever $g$ is, there are always non-trivial $g$-fixed vectors. This implies that $H^{1}(C,\mathbb{C})=0$. Hence $\sum_{p\in\mathbb{P}^{1}}(I_{\left<1\right>}-I_{\left<h_{p}\right>})-2I_{\left<1\right>}+2\mathbbm{1}=0$. But no combination will give this equality. Hence such $g$ does not exist. 
\end{proof}

{\bf Example 6 ($A_{5}$).} $G=A_{5}$ acts faithfully and symplectically on the K3 surface $\sum_{1}^{5}X_{i}=\sum_{1}^{6}X^{2}_{i}=\sum_{1}^{5}X^{3}_{i}=0$ in $\mathbb{P}^{5}$ via permutation action of the first 5 coordinates on $\mathbb{P}^{5}$. Then by Theorem \ref{curve1}, a $G$-stable curve $C$ with nodal singularities in an integral linear system can contribute to the representation $[e(\overline{J^{d}\mathcal{C}})]$ only if $\tilde{C}$ is rational.

\begin{prop}
Let $C$ be a smooth projective curve over $\mathbb{C}$ with a faithful $G=A_{5}$-action. If there exists $g\in A_{5}$ such that ${C}/\left<g\right>=\mathbb{P}^{1}$, then $C$ must be a smooth rational curve.
\end{prop}

\begin{proof} 
We have the following character table for $A_{5}$.
\begin{center}
\begin{tabular}{c||ccccc}
& $1A_{1}$ & $2B_{15}$ & $3A_{20}$ & $5A_{12}$ & $5B_{12}$ \\
\hline\hline
$\chi_{1}$ & 1&1&1&1&1\\
$\chi_{2}$ & $1^{(3)}$&$1(-1)^{(2)}$&$1\omega\bar{\omega}$ & $1\zeta_{5}\zeta_{5}^{4}$&$1\zeta_{5}^{2}\zeta_{5}^{3}$\\
$\chi_{3}$&$1^{(3)}$& $1(-1)^{(2)}$&$1\omega\bar{\omega}$ & $1\zeta_{5}^{2}\zeta_{5}^{3}$&$1\zeta_{5}\zeta_{5}^{4}$\\
$\chi_{4}$&$1^{(4)}$& $1^{(2)}(-1)^{(2)}$&$1^{(2)}\omega\bar{\omega}$&$\zeta_{5}...\zeta_{5}^{4}$&$\zeta_{5}...\zeta_{5}^{4}$\\
$\chi_{5}$&$1^{(5)}$& $1^{(3)}(-1)^{(2)}$&$1\omega^{(2)}\bar{\omega}^{(2)}$&$1\zeta_{5}...\zeta_{5}^{4}$&$1\zeta_{5}...\zeta_{5}^{4}$
\end{tabular}
\end{center}

For induced representations, we have
\[
\begin{aligned}
I_{2}&=\chi_{1}+\chi_{2}+\chi_{3}+2\chi_{4}+3\chi_{5}\\
I_{3}&=\chi_{1}+\chi_{2}+\chi_{3}+2\chi_{4}+\chi_{5}\\
I_{5}&=\chi_{1}+\chi_{2}+\chi_{3}+\chi_{5}
\end{aligned}
\]

Now suppose there exists $g\in G$ such that ${C}/\left<g\right>=\mathbb{P}^{1}$. Then we have $H^{1}(C,\mathbb{C})^{g}=H^{1}(C/\left<g\right>,\mathbb{C})=0$. But for $\chi_{1}$, $\chi_{2}$, $\chi_{3}$ and $\chi_{5}$, whatever $g$ is, there are always non-trivial $g$-fixed vectors. This implies that $H^{1}(C,\mathbb{C})$ does not contain $\chi_{1}$, $\chi_{2}$, $\chi_{3}$ and $\chi_{5}$ at all. Hence the coefficients of $\chi_{1}$, $\chi_{2}$, $\chi_{3}$ and $\chi_{5}$ in $\sum_{p\in\mathbb{P}^{1}}(I_{\left<1\right>}-I_{\left<h_{p}\right>})-2I_{\left<1\right>}+2\mathbbm{1}$ must be 0. This gives us only one possibility: $H^{1}(C,\mathbb{C})=I_{1}-I_{2}-I_{3}-I_{5}+2\mathbbm{1}=0$, which implies $C$ has genus 0. 
\end{proof}

{\bf Example 7 ($S_{5}$).} $G=S_{5}$ acts faithfully and symplectically on the K3 surface $\sum_{1}^{5}X_{i}=\sum_{1}^{6}X^{2}_{i}=\sum_{1}^{5}X^{3}_{i}=0$ in $\mathbb{P}^{5}$ via permutation action of the first 5 coordinates on $\mathbb{P}^{5}$. Then by Theorem \ref{curve1}, a $G$-stable curve $C$ with nodal singularities in an integral linear system can contribute to the representation $[e(\overline{J^{d}\mathcal{C}})]$ only if $\tilde{C}$ has genus 4. In this case $\tilde{C}$ has the largest possible automorphism group for a genus $4$ curve and $\tilde{C}$ is Bring's curve.

\begin{prop} 
Let $C$ be a smooth projective curve over $\mathbb{C}$ with a faithful $G=S_{5}$-action. If there exists $g\in S_{5}$ such that ${C}/\left<g\right>=\mathbb{P}^{1}$, then $C$ has genus $4$ and $g$ has order $5$. In particular, $C$ is Bring's curve.
\end{prop}

\begin{proof} 
We have the following character table for $S_{5}$.

\begin{center}
\begin{tabular}{c||cccc}
& $1A_{1}$ & $2A_{10}$ & $2A_{15}$ & $3A_{20}$\\
\hline\hline
$\chi_{1}$ & 1&1&1&1\\
$\chi_{2}$ & 1&-1&1&1\\
$\chi_{3}$&$1^{(4)}$& $1^{(3)}(-1)$&$1^{(2)}(-1)^{(2)}$&$1^{(2)}\omega\bar{\omega}$\\
$\chi_{4}$&$1^{(4)}$& $1(-1)^{(3)}$&$1^{(2)}(-1)^{(2)}$&$1^{(2)}\omega\bar{\omega}$\\
$\chi_{5}$&$1^{(5)}$& $1^{(3)}(-1)^{(2)}$&$1^{(3)}(-1)^{(2)}$&$1\omega^{(2)}\bar{\omega}^{(2)}$\\
$\chi_{6}$&$1^{(5)}$& $1^{(2)}(-1)^{(3)}$&$1^{(3)}(-1)^{(2)}$&$1\omega^{(2)}\bar{\omega}^{(2)}$\\
$\chi_{7}$&$1^{(6)}$& $1^{(3)}(-1)^{(3)}$&$1^{(2)}(-1)^{(4)}$&$1^{(2)}\omega^{(2)}\bar{\omega}^{(2)}$
\end{tabular}
\end{center}
\begin{center}
\begin{tabular}{c||ccc}
& $4A_{30}$ & $5A_{24}$ & $6A_{20}$\\
\hline\hline
$\chi_{1}$ & 1&1&1\\
$\chi_{2}$ &-1&1&-1\\
$\chi_{3}$ & $1{(-1)}i(-i)$&$\zeta_{5}...\zeta_{5}^{4}$&$1(-1)\omega\bar{\omega}$\\
$\chi_{4}$ & $1{(-1)}i(-i)$&$\zeta_{5}...\zeta_{5}^{4}$&$1(-1)(-\omega)(-\bar{\omega})$\\
$\chi_{5}$ & $1(-1)^{(2)}i(-i)$&$1\zeta_{5}...\zeta_{5}^{4}$&$1\omega\bar{\omega}(-\omega)(-\bar{\omega})$\\
$\chi_{6}$ & $1^{(2)}(-1)i(-i)$&$1\zeta_{5}...\zeta_{5}^{4}$&$-1\omega\bar{\omega}(-\omega)(-\bar{\omega})$\\
$\chi_{7}$ & $1(-1)i^{(2)}(-i)^{(2)}$&$1^{(2)}\zeta_{5}...\zeta_{5}^{4}$&$1(-1)\omega\bar{\omega}(-\omega)(-\bar{\omega})$
\end{tabular}
\end{center}

For induced representations, we have
\[
\begin{aligned}
I_{2A}&=\chi_{1}+3\chi_{3}+\chi_{4}+3\chi_{5}+2\chi_{6}+3\chi_{7}\\
I_{2B}&=\chi_{1}+\chi_{2}+2\chi_{3}+2\chi_{4}+3\chi_{5}+3\chi_{6}+2\chi_{7}\\
I_{3}&=\chi_{1}+\chi_{2}+2\chi_{3}+2\chi_{4}+\chi_{5}+\chi_{6}+2\chi_{7}\\
I_{4}&=\chi_{1}+\chi_{3}+\chi_{4}+\chi_{5}+2\chi_{6}+\chi_{7}\\
I_{5}&=\chi_{1}+\chi_{2}+\chi_{5}+\chi_{6}+2\chi_{7}\\
I_{6}&=\chi_{1}+\chi_{3}+\chi_{4}+\chi_{5}+\chi_{7}
\end{aligned}
\]

Now since there exists $g\in G$ such that ${C}/\left<g\right>=\mathbb{P}^{1}$, we have $H^{1}(C,\mathbb{C})^{g}=H^{1}(C/\left<g\right>,\mathbb{C})=0$. But for $\chi_{1}$, $\chi_{5}$ and $\chi_{7}$, whatever $g$ is, there are always non-trivial $g$-fixed vectors. This implies that $H^{1}(C,\mathbb{C})$ does not contain $\chi_{1}$, $\chi_{5}$ and $\chi_{7}$ at all. Hence the coefficients of $\chi_{1}$, $\chi_{5}$ and $\chi_{7}$ in $\sum_{p\in\mathbb{P}^{1}}(I_{\left<1\right>}-I_{\left<h_{p}\right>})-2I_{\left<1\right>}+2\mathbbm{1}$ must be 0. This gives us only two possibilities: $H^{1}(C,\mathbb{C})=I_{1}-I_{2}-I_{4}-I_{5}+2\mathbbm{1}=2\chi_{4}$ or $H^{1}(C,\mathbb{C})=I_{1}-I_{2}-I_{5}-I_{6}+2\mathbbm{1}=2\chi_{4}+2\chi_{6}$. For the second case, we notice that whatever conjugacy class $g$ belongs to, there always exists $g$-fixed vectors in $2\chi_{4}+2\chi_{6}$. Hence $H^{1}(C,\mathbb{C})=I_{1}-I_{2}-I_{4}-I_{5}+2\mathbbm{1}=2\chi_{4}$. It follows that $C$ has genus $4$ and $g$ has order $5$. 
\end{proof}

\begin{center}
\sc{Acknowledgements}
\end{center}
\vspace{0.1 in}
I thank my advisor Professor Michael Larsen for his guidance and valuable discussions throughout this work, and in particular, for suggesting the problem. I also thank Professor Jim Bryan for several helpful suggestions and comments. Finally, I thank both referees for their careful reading of this manuscript and many helpful suggestions. See \cite{Z21} for the version of record.


\begin{thebibliography}{DD2}

\bibitem[Ale04]{A04}Alexeev, V.: Compactified Jacobians and Torelli map. \textit{Publ.\ RIMS,\ Kyoto\ Univ.} \textbf{40} (2004), 1241--1265.

\bibitem[AK76]{AK76}Altman, A.; Kleiman, S.: Compactifying the Jacobian, \textit{Bull.\ Amer.\ Math.\ Soc.} \textbf{82(6)} (1976), 947--949.

\bibitem[AKMW02]{AKMW02} Abramovich, D.; Karu, K.; Matsuki, K.; W\l odarczyk, J.: Torification and factorization of birational maps. \textit{J.\ Amer.\ Math.\ Soc.} \textbf{15(3)} (2002), 531--572.

\bibitem[AST11]{AST11}Artebani, M.; Sarti, A.; Taki, S.: K3 surfaces with non-symplectic automorphisms of prime order. \textit{Math.\ Z.} \textbf{268} (2011), 507--533.

\bibitem[Bea99]{Bea99}Beauville, A.: Counting rational curves on K3 surfaces. \textit{Duke\ Math.\ J.} \textbf{97} (1999), 99--108.

\bibitem[Bli11]{Bl11}Blickle, M.: A short course on geometric motivic integration. Motivic integration and its interactions with model theory and non-Archimedean geometry. Volume I, London Math. Soc. Lecture Note Ser., 383, Cambridge Univ. Press, Cambridge, 2011, 189--243.

\bibitem[BG19]{BG19}Bryan, J.; Gyenge, \'A.: G-fixed Hilbert schemes on K3 surfaces, modular forms, and eta products. arXiv:1907.01535. 

\bibitem[BLR90]{BLR90}Bosch, Siegfried; L\"utkebohmert, Werner; Raynaud, Michel: N\'eron models, Ergebnisse der Mathematik und ihrer Grenzgebiete (3)[Results in Mathematics and Related Areas(3)], vol. 21, pp. x+325, Springer-Verlag, Berlin, 1990.

\bibitem[BO18]{BO18}Bryan, J.; Oberdieck, G.: CHL Calabi-Yau threefolds: Curve counting, Mathieu moonshine and Siegel modular forms. arXiv:1811.06102.

\bibitem[Car85]{Ca85}Carter, R.W.: Finite Groups of Lie Type: Conjugacy Classes and Complex Characters, Wiley, 1985.

\bibitem[CS17]{CS17}Cohen, H.; Str\"omberg, F: Modular Forms: A Classical Approach, Graduate Studies in Math. 179, Amer. Math. Soc., 2017.

\bibitem[DH88]{DH88}Diaz, Steven; Harris, Joe: Ideals associated to deformations of singular plane curves.
\textit{Trans.\ Amer.\ Math.\ Soc.} \textbf{309} (1988), no.\ 2, 433--468.

\bibitem[DL76]{DL76}Deligne, P.; Lusztig, G.: Representations of reductive groups over finite fields. \textit{Ann. of Math.} \textbf{103} (1976), 103--161.

\bibitem[DL99]{DL99}Denef, J.; Loeser, F.: Germs of arcs on singular algebraic varieties and motivic integration.
\textit{Invent. Math.} \textbf{135} (1999), no. 1, 201–232.

\bibitem[DL02]{DL02}Denef, J.; Loeser, F.: Lefschetz numbers of iterates of the monodromy and truncated arcs. \textit{Topology} \textbf{41} (2002), no. 5, 1031–1040. 

\bibitem[EG00]{EG00}Ellingsrud, G.; G\"ottsche, L.: Hilbert schemes of points and Heisenberg algebras. School on Algebraic Geometry (Trieste, 1999), 59–100, ICTP Lect. Notes, 1, Abdus Salam Int. Cent. Theoret. Phys., Trieste, 2000. 

\bibitem[EGK00]{EGK00}Esteves, E.; Gagn\'e, M.; Kleiman, S.: Abel maps and presentation schemes. \textit{Comm. in Algebra} \textbf{28(12)} (2000), 5961--5992.

\bibitem[G\"ot90]{Go90}G\"ottsche, L.: The Betti numbers of the Hilbert scheme of points on a smooth projective surface. \textit{Math.\ Ann.} \textbf{286} (1990), 193--207.

\bibitem[G\"ot94]{Go94}G\"ottsche, L.: Hilbert Schemes of Zero-Dimensional Subschemes of Smooth Varieties, LNM 1572, Springer-Verlag Berlin Heidelberg, 1994.

\bibitem[Gro96]{Gr96}Grojnowski, I.: Instantons and affine algebras. I. The Hilbert scheme and vertex operators, \textit{Math.\ Res.\ Lett.} \textbf{3(2)} (1996), 275--291.

\bibitem[GS93]{GS93}G\"ottsche, L.; Soergel, W.: Perverse sheaves and the cohomology of Hilbert schemes of smooth algebraic surfaces. \textit{Math.\ Ann.} \textbf{296} (1993), 235--245.

\bibitem[GS07]{GS07} van Geemen, B.; Sarti, A.: Nikulin involutions on K3 surfaces. \textit{Math.\ Z.} \textbf{255} (2007), 731--753.

\bibitem[Huy16]{H16}Huybrechts, D.: Lectures on K3 Surfaces, Cambridge studies in advanced mathematics 158, Cambridge University Press, 2016.

\bibitem[Iar77]{I77}Iarrobino, A.A.: Punctual Hilbert Schemes. Mem. Amer. Math. Soc. \textbf{10} , no. 188, viii+112 pp, 1977.

\bibitem[Loe09]{L09}Loeser, F.: Seattle lectures on motivic integration. Algebraic geometry—Seattle 2005. Part 2, Proc. Sympos. Pure Math., 80, Part 2, Amer. Math. Soc., Providence, RI, 2009, 745--784.

\bibitem[LN20]{LN20}L$\hat{e}$, Q.; Nguyen, H.: Equivariant motivic integration and proof of the integral identity conjecture for regular functions. \textit{Math. Ann.} \textbf{376(3-4)} (2020), 1195–1223.

\bibitem[Mat00]{Ma00}Matsuki, K.: Lectures on factorization of birational maps. arXiv:math/0002084.

\bibitem[Muk88]{Mu88}Mukai, S.: Finite groups of automorphisms of K3 surfaces and the Mathieu group. \textit{Invent.\ math.} \textbf{94} (1988), 183--221.

\bibitem[Nak97]{N97}Nakajima, H.: Heisenberg algebra and Hilbert schemes of points on projective surfaces. \textit{Ann.\ of\ Math.} (2) \textbf{145(2)} (1997), 379--388.

\bibitem[Nak99]{N99}Nakajima, H.: Lectures on Hilbert Schemes of Points on Surfaces. University Lecture Series, 18. American Mathematical Society, Providence, RI, 1999. xii+132 pp.

\bibitem[PS08]{PS08}Peters, C.; Steenbrink, J.: Mixed Hodge structures. Results in Mathematics and Related Areas. 3rd Series. A Series of Modern Surveys in Mathematics, 52. Springer-Verlag, Berlin, 2008.

\bibitem[Ser79]{Se79}Serre, J.P.: Local Fields, GTM 67, Springer-Verlag New York, 1979.

\bibitem[Ser88]{Se88}Serre, J.P.: Algebraic Groups and Class Fields, GTM 117, Springer-Verlag New York, 1988.

\bibitem[YZ96]{YZ96}Yau, S.T.; Zaslow, E.: BPS states, String duality, and Nodal curves on K3. \textit{Nuclear\ Physics\ B} \textbf{471(3)} (1996), 503--512.

\bibitem[Zha21]{Z21}Zhan, S.: Counting rational curves on K3 surfaces with finite group actions. \textit{Int. Math. Res. Not. IMRN} (2021), https://doi.org/10.1093/imrn/rnaa320.
\end{thebibliography}
\end{document}